\documentclass[12pt]{article}

\usepackage[dvips]{graphicx}

\usepackage{amssymb,amsthm,amsmath,latexsym}

\newcommand{\od}{\overset{\mbox{\rm def}}=}

\theoremstyle{definition}
\newtheorem{definition}{Definition}

\theoremstyle{plain}
\newtheorem{lemma}{Lemma}

\theoremstyle{plain}
\newtheorem{theorem}{Theorem}

\theoremstyle{plain}
\newtheorem{corollary}{Corollary}

\theoremstyle{definition}
\newtheorem{remark}{Remark}

\theoremstyle{definition}
\newtheorem{example}{Example}

\begin{document}

\title{Confluence of the nonlinear waves in the Stefan problem
with undercooling}

\author{V.G.~Danilov\thanks{Moscow Technical University of
Communication and Informatics,\hfill\break e-mail:
danilov@miem.edu.ru }\and V.Yu.~Rudnev\thanks{Moscow Technical
University of Communication and Informatics, e-mail:
pm@miem.edu.ru \hfill\break This work was supported by the Russian
Foundation for Basic Research under grant 05-01-00912.}}

\date{}
\maketitle

\abstract{We assume that the Stefan problem with undercooling has a
classical solution until the moment of contact of free boundaries
and the free boundaries have continuous velocities until the moment
of contact. Under these assumptions, we construct a smooth
approximation of the global solution of the Stefan problem with
undercooling, which, until the contact, gives the classical solution
mentioned above and, after the contact, gives a solution which is
the solution of the heat equation.}

\section{Introduction}\label{S:intro}

In this paper, we study the confluence of free boundaries in the
Stefan problem with kinetic undercooling in the spherical
symmetric case.

In the preset paper, we study the problem in the domain
$Q=\Omega\times[0,t_{1}]$, where $\Omega=[R_{1},R_{2}]$ is an
spherical layer in the spherical coordinates. We assume that the
domain $\Omega$ is divided into three layers $\Omega^{+}_{1,2}(t)$
and $\Omega^{-}(t)$ as follows:
\begin{equation*}
\Omega^{+}_{1}(t)=[R_{1},r_{1}(t)], \quad
\Omega^{-}(t)=[r_{1}(t),r_{2}(t)], \quad
\Omega^{+}_{2}(t)=[r_{2}(t), R_{2}],
\end{equation*}
where $r_{i}(t)$, $i=1,2$, are the free boundaries of phases "$+$"
and "$-$". We assume that the phase~"$+$" occupies the layers
$\Omega^{+}_{1,2}(t)$ and the phase~"$-$" occupies the layer
$\Omega^{-}(t)$.

We shall construct a smooth approximation of solutions
of the Stefan problem with undercooling
under the assumption that the motion of the free boundary
is the motion of the front of a nonlinear wave
and the confluence of free boundaries
is interaction of solitary nonlinear waves.
The possibility of this interpretation is given  by
the models of phase field~\cite{12} proposed by G.~Caginalp.
In fact, the choice of the method for approximating
the limit Stefan problems with undercooling is unessential for
us, because we do not prove that the approximations thus
constructed are close to the corresponding solutions
of the phase field system.

For example, we could use the definition of the generalized
solution of the limit problems including the order function
(the nonlinear wave)~\cite{111}.

Here our considerations are based on the following simple fact.
Suppose that there are two families of solutions
(exact and approximate solutions)
of some problem depending on a small parameter~$\varepsilon$.
Suppose that both these families have the properties that
permit passing to the limit as $\varepsilon\to 0$
in the weak sense.
Suppose also that the family of approximate solutions
satisfies a problem with a right-hand side
small as $\varepsilon\to 0$ in the weak sense.
Then the weak limits of both families are solutions
of the same limit problem and
if the latter has a unique solution, then the difference
between the exact and approximate solutions tends to~$0$
as $\varepsilon\to0$ at least in the weak sense.

We recall that the phase field system has the form
\begin{equation}\label{PF:1}
L\theta= -\frac{\partial u}{\partial t},\qquad
\varepsilon Lu-\frac{u-u^{3}}{\varepsilon}-\theta=0,
\end{equation}
where
\begin{equation*}\label{L:1}
L=\frac{\partial}{\partial
t}-\frac{1}{r^{2}}\frac{\partial}{\partial
r}\left(r^{2}\frac{\partial}{\partial r}\right), \qquad r\in
[R_{1},R_{2}], \quad t\in[0,t_{1}].
\end{equation*}

The function $\theta=\theta(r,t,\varepsilon)$ has the meaning of
the temperature, and the function $u=u(r,t,\varepsilon)$, which is
called the order function, determined the phase state of the
medium: $u\simeq-1$ corresponds to the phase~"$-$" in the layer
$\Omega^{-}(t)$, and $u\simeq 1$ corresponds to the phase~"$+$" in
the layers  $\Omega^{-}_{1,2}(t)$.

Passing to the limit as $\varepsilon\to 0$ in (\ref{PF:1}), we
obtain the Stefan problem with undercooling.

This passage to the limit is possible, for example, in the case
where the corresponding limit problems have classical solutions.
In this case, the weak limits as $\varepsilon\to 0$ of solutions
(\ref{PF:1}) give these solutions~\cite{12,2,3}.

We denote
\begin{equation}\label{sigma}
\theta(r,t)=\frac{\sigma(r,t)}{r}.
\end{equation}
Then the phase field system (\ref{L:1}) we can write in the form
\begin{align}
&\frac{\partial\sigma}{\partial
t}-\frac{\partial^{2}\sigma}{\partial r^{2}}=
-r\frac{\partial u}{\partial t},\label{PhHeat}\\
&\varepsilon
Lu=\frac{u-u^{3}}{\varepsilon}+\frac{\sigma}{r}.\label{PhAC}
\end{align}

The smooth approximations of solutions of the Stefan problem with
kinetic undercooling constructed in this paper are approximate (in
the above sense) solutions of the phase field system and they
admit a weak passage to the limit as $\varepsilon\to 0$. In this
case, we obtain the limit problems (and their solutions)
describing the process of confluence of free boundaries. To
construct these approximations, we use the assumption that the
classical sharp fronted solution of the Stefan problem with
kinetic undercooling exists until the confluence of free
boundaries begins. This is a natural assumption in our paper
(otherwise, it is not clear the confluence of what is considered)
and can be proved\footnote{A.~M.~Meirmanov. A private
communication.}.

We repeat that, under this assumption, we construct an
approximation of solution of the limit problem, which is smooth
for $\varepsilon>0$ and uniformly bounded in $\varepsilon\geqslant
0$. As is well known, in a similar problem about the propagation
of shock waves, the existence of such an approximation
distinguishes a unique solution.

By $t^{*}\in[0,t_{1})$ we denote the instant of confluence of free
boundaries. Then, for any $t\leqslant t^{*}-\delta$ for all
$\delta>0$, we see that the limit problems have solutions. The
asymptotic solution of system (\ref{PhHeat}),(\ref{PhAC}) has the
form
\begin{align}\label{AS:1}
\sigma^{as}_{\varepsilon}&=\bar{\sigma}^{-}(r,t)+
\left(\bar{\sigma}^{+}(r,t)-\bar{\sigma}^{-}(r,t)\right)
\omega_{1}\left(\frac{\hat{r}(t)-r}{\varepsilon}\right)
\omega_{1}\left(\frac{r+\hat{r}(t)}{\varepsilon}\right),
\\
\label{AS:2} u^{as}_{\varepsilon}&=
1+\omega_{0}\left(\frac{-r-\hat{r}(t)}{\varepsilon}\right)
\\
&\qquad +\omega_{0}\left(\frac{r-\hat{r}(t)}{\varepsilon}\right)
+\varepsilon\left[\frac{\sigma^{as}_{\varepsilon}}{2}+
\omega\left(t,\frac{r+\hat{r}(t)}{\varepsilon},
\frac{r-\hat{r}(t)}{\varepsilon}\right)\right].\notag
\end{align}
Here $\omega_{1}(z)\to 0,1$ as $z\to\mp\infty$,
$\omega^{(k)}_{1}(z)\in\mathbb{S}(\mathbb{R}^{1}_{z})$ for $k>0$,
$\hat{\varphi}(t)$ is a smooth function, $\omega_{0}(z)=\tanh(z)$,
and $\omega(t,z_{1},z_{2})\in
C^{\infty}([0,t^{*}];\,\mathbb{S}(\mathbb{R}^{2}_{z}))$. By
$\mathbb{S}(\mathbb{R}^n)$ we denote the Schwartz space of smooth
rapidly decreasing functions. If the initial data for (\ref{PF:1})
has the form (\ref{AS:1}), (\ref{AS:2}) at $t=0$, then, for
$t\leqslant t^{*}-\delta$, we have the estimate
\begin{equation*}
\|u-u^{as}_{\varepsilon};C(0,T;L^{2}(\mathbb{R}^{1})\|+\|\sigma-
\sigma^{as}_{\varepsilon};\mathcal{L}^{2}(Q)\|\leqslant
c\varepsilon^{\mu},\qquad \mu\geqslant  3/2,
\end{equation*}
where $(\sigma,u)$ is a solution of system
(\ref{PhHeat}),(\ref{PhAC}) (see~\cite{1,2}). Here
$Q=\Omega\times[0,t^{*}-\delta)$, and the constant $c$ is
independent of~$\varepsilon$.

The main obstacle to the construction of solutions of the form
(\ref{AS:1}), (\ref{AS:2}), which could be used to describe the
confluence of free boundaries, is the fact that, instead of an
ordinary differential equation whose solution is the function
$\omega_{0}(z)$~\cite{12,2}, in the case of confluence of free
boundaries, we must deal with a partial differential equation for
which the explicit form of the exact solution is unknown.

In the present paper, we use the technique of the weak asymptotics
method \cite{6,11}, which allows us to avoid this problem. Let us
explain several basic points.

\begin{definition}\label{WkMDef}
A family of functions $f(r,t,\varepsilon)\in L^1(Q)$ integrable
with respect to~$r$ for all $t\in[0,t_{1}]$ and for
$\varepsilon>0$ admits the estimate
$\mathcal{O}_{\mathcal{D}'}(\varepsilon^\nu)$ if, for any test
function $\zeta(r)\in C^{\infty}_{0}$, we have the estimate
\begin{equation}\label{WkMDef:1}
\left|\int_{\Omega}f(r,t,\varepsilon)\zeta(r)dx\right|\leqslant
C_{t_{1},\zeta}\varepsilon^{\nu},
\end{equation}
where the constant $C_{t_{1},\zeta}$ depends on~$t_{1}$ and the
test function $\zeta(r)$.
\end{definition}

Generalizing (\ref{WkMDef:1}), we shall say that the family of
distributions $f(r,t,\varepsilon)$ depending on~$t$
and~$\varepsilon$ as on parameters admits the estimate
$\mathcal{O}_{\mathcal{D}'}(\varepsilon^\nu)$ if, for any test
function $\zeta(r)$, we have the estimate
\begin{equation*}
\langle f(r,t,\varepsilon), \zeta(x)\rangle
=\mathcal{O}(\varepsilon^\nu),\quad 0\leqslant t \leqslant t_{1}.
\end{equation*}

\begin{example}
Let $\omega(z)\in \mathbb{S}(\mathbb{R}^{1})$,
$\int_{\mathbb{R}^{1}}\omega(z)\,dz=1$, and let
$x_{0}\in[l_{1},l_{2}]$. Then we have
\begin{equation*}
\omega\left(\frac{x}{\varepsilon}\right)-\delta(x)
=\mathcal{O}_{\mathcal{D}'}(\varepsilon), \qquad
\omega\left(\frac{x-x_{0}}{\varepsilon}\right)
=\delta(x-x_{0})\int_{\mathbb{R}^{1}}\omega(z)\,dz
+\mathcal{O}_{\mathcal{D}'}(\varepsilon).
\end{equation*}
\end{example}

\begin{example}
Suppose that $\omega_{i}(z) \in C^{\infty}$, \
$\left(\omega_{i}\right)'\in\mathbb{S}(\mathbb{R}^{1})$, \
$\lim_{z\to-\infty}\omega_{i}(z)=0$, \ and
$\lim_{z\to+\infty}\omega_{i}(z)=1$, $i=1,2$. Then we have
\begin{description}
\item{(a)}
\begin{equation*}
\omega_{i}\left(\frac{x-x_{i}}{\varepsilon}\right)-H(x-x_{i})
=\mathcal{O}_{\mathcal{D}'}(\varepsilon).
\end{equation*}
\item{(b)}
\begin{align*}
&\omega_{1}\left(\frac{x-x_{1}}{\varepsilon}\right)
\omega_{2}\left(\frac{x-x_{2}}{\varepsilon}\right)\\
&\qquad =B_{1}\left(\frac{\triangle
x}{\varepsilon}\right)H(x-x_{1}) +B_{2}\left(\frac{\triangle
x}{\varepsilon}\right)H(x-x_{2})
+\mathcal{O}_{\mathcal{D}'}(\varepsilon),
\end{align*}
\end{description}
where $\triangle x=x_{1}-x_{2}$, \ $B_{i}(\rho)\in C^{\infty}$, \
$B_{i}'(\rho)\in\mathbb{S}(\mathbb{R}^{1})$, \
$B_{1}(\rho)+B_{1}(\rho)=1$, \ $B_{1}(+\infty)=1$, and
$B_{1}(-\infty)=0$.
\end{example}

\begin{example}
The preceding relations easily imply the formula
\begin{align*}
H(x-x_1)H(x-x_2) =B\left(\frac{\triangle
x}{\varepsilon}\right)H(x-x_1)
+(1-B)H(x-x_2)+\mathcal{O}_{\mathcal{D}'}(\varepsilon),
\end{align*}
where $B(\rho)\in C^\infty$ is the function $B_1$ from the
preceding example.
\end{example}

\begin{example}
The following relation is a corollary of Definition~1:
$$
\bigg(\frac{d}{dx}\bigg)^m
\mathcal{O}_{\mathcal{D}'}(\varepsilon^\alpha)
=\mathcal{O}_{\mathcal{D}'}(\varepsilon^\alpha)
$$
for all $m\in\mathbb{Z}_+$ and $\alpha>0$. Generally speaking,
this is not true for the derivatives w.r.t.~$t$.
\end{example}

The relations in Examples~2\,(a), 2\,(b) and~4 are obvious, but,
for the reader's convenience, we prove some of these formulas in
Section~5.

Here we only note that, in view of items~(a) and~(b), one can say
that the product of approximations of the Heaviside functions or
of the Heaviside functions themselves is a linear combination with
accuracy up to small terms
($\sim\mathcal{O}_{\mathcal{D}'}(\varepsilon)$). The general
nonlinear functions of linear combinations of approximations of
the Heaviside functions can be linearized similarly. Precisely
this property underlies the constructive study of the interaction
between nonlinear waves with localized fast variations.

\begin{definition}\label{Def:1}
A pair of smooth functions $(\check{u}, \check{\sigma})$ is a weak
asymptotic solution of the phase field system
(\ref{PhHeat}),(\ref{PhAC}) if, for any test functions
$\zeta(r),\xi(r)\in C^{\infty}_{0}(\Omega)$, the following
relations hold:
\begin{align}\label{Defin:1}
&\int_{\Omega}\left(r\check{u}_{t}+\check{\sigma}_{t}\right)r^{2}\zeta
dr +\int_{\Omega}\check{\sigma}_{r}(r^{2}\zeta)_{r}dr =
\mathcal{O}(\varepsilon),
\\
\label{Defin:2} &\varepsilon\int_{\Omega}
\check{u}_{r}\check{u}_{t}\xi dr
+2\varepsilon\int_{\Omega}r\check{u}^{2}_{r} \xi
dr-\frac{\varepsilon}{2}\int_{\Omega}\check{u}^{2}_{r}(r^{2}
\xi)_{r} dr \\
&-\frac{1}{\varepsilon}\int_{\Omega} \left(\frac{\check{u}^{4}}{4}
-\frac{\check{u}^{2}}{2}+\frac{1}{4}\right)(r^{2}\xi)_{r} dr
+\int_{\Omega}\check{u}(r\check{\sigma}\xi)_{r}dr =
\mathcal{O}(\varepsilon^\mu), \qquad \mu\in(0,1/2).\notag
\end{align}
\end{definition}

The left-hand side of Eq.~(\ref{Defin:2}) is obtained by
multiplying equation (\ref{PhAC}) by $\check{u}_{r}$ and
integrating by parts. The reminders
$\mathcal{O}(\varepsilon^\alpha)$, $\alpha=1,\mu$, in the
right-hand sides of (\ref{Defin:1}) and (\ref{Defin:2}) must be
bounded locally in $t$, i.e., for $t\in[0,t_{1}]$, we have
\begin{equation*}
\max_{0\leqslant t\leqslant
t_{1}}\left|\mathcal{O}(\varepsilon)\right|\leqslant
C_{t_{1}}\varepsilon,\qquad C_{t_{1}}=\mathrm{const}.
\end{equation*}
This construction was introduced and analyzed in~\cite{4}.

The fact that a noninteger exponent appears in the right-hand side
of~(\ref{Defin:2}) is not directly related to the technique of the
weak asymptotics method, see Examples 1--4. The source of this
noninteger exponent is the nonsmoothness of the function $\theta$
(the temperature), which appears at the instant of confluence of
the free boundaries.

The paper is organized as follows.

In Section~2, we explain the structure of ansatzes of
approximations of the temperature and the order function. These
ansatzes are constructed under the assumption the classical sharp
fronted solution of the Stefan problem exists. At the beginning of
Section~2, we formulate what we exactly need.

Next, in Section~3, we substitute the constructed ansatzes in
system~(\ref{PhHeat}),(\ref{PhAC}) and derive equations for the
unknown functions contained in the ansatzes.

The results of these calculations are summarized in Theorems~1--3
in Section~3. The final result of this paper is formulated in
Theorem~4 in Section~3. It states that the assumption on the
existence of the classical solution is sufficient for constructing
formulas for the weak asymptotic solution of
system~(\ref{PhHeat}),(\ref{PhAC}) in the sense of
Definition~\ref{WkMDef}. Next, in Section~4, we analyze the
constructed formulas and derive the following effects:
\begin{description}
\item{(a)}
the weak asymptotic solution is smooth for $t>t^*$, the absolute
values of the free boundaries velocities are equal to each other
at the contact moment;
\item{(b)}
the temperature has a negative jump at the instant and at the point of confluence of the free boundaries, and this jump is equal to
$$
[\bar{\theta}]|_{\substack{r=r^*,\\t=t^*}} =-\frac{r_{10}(t^*)}4\lim_{t\to t^*-0}(r_{10t}^{2}+r_{20t}^{2}).
$$
where $r=r_{i0}(t)$ other positions of the free boundaries till the interaction.
\end{description}

In particular, it follows from (a) and (b) that the velocities of
the free boundaries have jumps at the point of contact.

These effects can also be discovered in the numerical analysis of
the process of confluence of free boundaries. Some helpful
technical results are given in Section~5. We note that the only
example known to the authors, where the confluence of the free
boundaries is studied, is given in \cite{13}.

In general, this paper turned out to be technically rather
complicated and long in spite of the fact that the details in
several justifications in Section~5 were omitted.

Although, as was shown above, the results obtained here are a
necessary step in the study of the multidimensional problem, this
paper shows that the following classification can be introduced:

(1) problems in which the confluence of free boundaries leads to
disappearance of one of the phases;

(2) problems in which the domain occupied by one of the phases
changes its connectivity, but the number of phases remains the
same.

In this paper, we consider an example precisely from the first
class of problems. As was mentioned above, we do not justify the
asymptotics of the constructed solution. But this can be done
based on our constructions. It is easy to see that the main role
here is played by the (not proved) existence of the classical
solution up to the moment of confluence of the free boundaries.
Our construction of the existnece under this assumption reduced
the justification to estimating the soluton of the heat equation
with the right-hand side $f_\varepsilon$ admitting the estimate
$$
f_\varepsilon=O_{\mathcal{D}'}(\varepsilon^\mu),\qquad
\mu\in(0,1/2)
$$
and with zero initial and boundary conditions. An analysis of the
structure of this right-hand side shows that
$f_\varepsilon\cdot\varepsilon^{-\mu}$ as $\varepsilon\to0$ is a
linear combination of functions $\delta'(r-r_{i})$ and
$\delta(r-r_{i})$, $i=1,2$ with coefficients depending on
$t,\tau$, and these coefficients converge fast to zero as
$\tau\to\pm\infty$.

Hence we can conclude that for $r\ne r_{i}$ the solution of this
heat equation belongs to $C^\infty$ for $\varepsilon\geq0$ and
admits the estimate $\mathcal{O}(\varepsilon^\mu)$. In the whole
domain $\Omega\times[0,T]$, a rough analysis based on general
theorems~\cite{14} shows that the solution belongs to
$W^{-\delta}_2$, $\delta>0$, and has the estimate
$\mathcal{O}(\varepsilon^\mu)$ in the norm of this space.  In
particular, the weak limit of the constructed weak asymptotic
solution is equal to the exact global solution of the heat
equation in the phase field system.

\section{Ansatz of the approximation of the solution
of the Stefan problem with undercooling}\label{S:Constr}

In our case, the Stefan problem with kinetic undercooling has the
form
\begin{align}
&\frac{\partial\overline{\sigma}}{\partial t}
=\frac{\partial^{2}\overline{\sigma}}{\partial r^{2}},\qquad
r\in[R_{1},R_{2}],\quad r\neq \hat{r}_{i}(t),\quad
i=1,2,\label{SH:1}\\
&\left.\frac{\overline{\sigma}}{r}\right|_{r=\hat{r}_{i}(t)}=(-1)^{i+1}\left(
\varkappa_{1}\hat{r}_{i}'(t)+\frac{\varkappa_{2}}{\hat{r}_{i}(t)}\right),\label{SGT:1}\\
&\left.\left[\frac{\partial\overline{\sigma}}{\partial r}\right]
\right|_{r=\hat{r}_{i}(t)}=(-1)^{i+1}\hat{r}_{i}^{3}(t)\hat{r}'_{i}(t).\label{SSt:1}
\end{align}
Relations (\ref{SH:1})--(\ref{SSt:1}) are supplemented with
(Dirichlet or Neumann) boundary conditions for $r=R_{i}$, $i=1,2$,
and with a consistent initial condition.

Obviously, problem (\ref{SH:1})--(\ref{SSt:1}) can be written as
\begin{align}\label{SH:R}
&\frac{\partial\overline{\sigma}}{\partial t}
-\frac{\partial^{2}\overline{\sigma}}{\partial
r^{2}}=\hat{r}^{3}_{2}\hat{r}_{2t}\delta(r-\hat{r}_{2})-
\hat{r}^{3}_{1}\hat{r}_{1t}\delta(r-\hat{r}_{1}),
\\
&\left.\frac{\overline{\sigma}}{r}\right|_{r=\hat{r}_{i}(t)}=(-1)^{i+1}
\left(\varkappa_{1}\hat{r}_{i}'(t)+\frac{\varkappa_{2}}{\hat{r}_{i}(t)}\right).\notag
\end{align}
We shall assume that the initial conditions are chosen so that the
problem in question has the classical solution,
$\hat{r}_{2}(0)-\hat{r}_{1}(0)>0$, and there exists a
$t=t^{*}\in[0,t_{1}]$ such that
$\hat{r}_{1}(t^{*}-0)=\hat{r}_{2}(t^{*}-0)$. More precisely, we
assume that
\begin{description}
\item{(i)}
the limits $\lim_{t\to t^{*}-0}\hat{r}_{it}$, $i=1,2$ exist;
\item{(ii)}
\begin{align*}
\overline{\sigma}(r,t)&\in C^{1}\left([0,t^{*}),
C^{2}\left(\mathring{\Omega}{}^{-}(t)\cup
\mathring{\Omega}{}^{+}_{1}(t)
\cup\mathring{\Omega}{}^{+}_{2}(t)\right)\right)\\
&\cap\left(C^{1}\left(\Omega^{-}(t)\right)\cup
C^{1}\left(\Omega^{+}_{1}(t)\right)\cup
C^{1}\left(\Omega^{+}_{2}(t)\right)\right)\cap C(\Omega)
\end{align*}
(here $\mathring{D}$ denotes the interior of the domain~$D$),
cf.~\cite{9}.
\end{description}

As will be shown below, these assumptions permit describing the
interaction (the confluence of free boundaries) constructively.

We shall construct this approximation (ansatz) as a weak
asymptotic solution of system~(\ref{PhHeat}),(\ref{PhAC}).

First, we introduce the ansatz of the order function in
(\ref{SH:1})--(\ref{SSt:1}). It has the form
\begin{align}\label{AU:1}
 \check{u}&=\frac{1}{2}\left[1+
 \omega_{0}\left(\beta\frac{r_{1}-r}{\varepsilon}\right)+
 \omega_{0}\left(\beta\frac{r-r_{2}}{\varepsilon}\right)\right.\\
          &   -\left.
  \omega_{0}\left(\beta\frac{r_{1}-r}{\varepsilon}\right)
  \omega_{0}\left(\beta\frac{r-r_{2}}{\varepsilon}\right)\right]\notag.
\end{align}

Here the unknowns are the functions $\beta$ and
$r_{i}=r_{i}(t,\varepsilon)$, $i=1,2$. In what follows, we write
them more precisely. Now we only note that if $\beta>0$, then, for
$r_{1}<r_{2}$, the function $\check{u}$ coincides up to
$\mathcal{O}(\varepsilon^{N})$ with the sum of the first three
terms in the right-hand side of (\ref{AU:1}), for $r_{1}=r_{2}$,
we have $\check{u}=1+\mathcal{O}_{\mathcal{D}'}(\varepsilon)$, and
for $r_{1}>r_{2}$, we have
$\check{u}=1+\mathcal{O}(\varepsilon^N)$. These relations can
easily be verified directly. We note that, under the above
assumptions, we can continue the functions $\hat{r}_{i}(t)$ to the
interval $[0,t_{1}]$ preserving the smoothness and the sign of the
derivatives. We choose such continuations and denote them by
$r_{i0}(t)$. Now, we can write the functions
$r_{i}(t,\varepsilon)$ and $\beta(t,\varepsilon)$ more precisely.
Namely, we set
\begin{equation}\label{PHI:1}
r_{i}(t,\varepsilon) =r_{i0}(t)+\psi_{0}(t)r_{i1}(\tau,t),
\end{equation}
where
\begin{equation}\label{PhiTau:1}
\psi_{0}(t)=r_{20}(t)-r_{10}(t),\quad \tau=\psi_{0}/\varepsilon.
\end{equation}
We note that
\begin{align*}
&\tau\to \infty,\quad t<t^{*},\\
&\tau\to -\infty,\quad t>t^{*},
\end{align*}
as $\varepsilon\to 0$. Thus, the value of the variable~$\tau$
characterizes the process of confluence of free boundaries. Also,
we introduce the functions
\begin{equation}\label{R1PH1}
\psi(t,\varepsilon)=r_{2}(t,\varepsilon)-r_{1}(t,\varepsilon),
\quad \rho(\tau)=\psi(t,\varepsilon)/\varepsilon
\end{equation}

We set
\begin{equation}\label{BETADef:1}
\beta(\tau)=1+\beta_{1}(\tau).
\end{equation}
In addition, we assume that
\begin{align}
&\beta^{(\alpha)}_{1}(\tau),\ \varphi_{i1}^{(\alpha)}(\tau)\to
\mathcal{O}(\tau^{-1}),\quad \tau\to\infty, \quad
\alpha=0,1.\label{BETAPHI:1}\\
&0<\delta_{1}<\beta_{1}(\tau)<\delta_{2},\quad \delta_{1},\
\delta_{2}=\mathrm{const},\notag\\ &r_{i1}\to r^{-}_{i1},\quad
r_{i1}-r^{-}_{i1} =\mathcal{O}(|\tau|^{-1}),\quad\tau\to
-\infty,\quad r^{-}_{i1}=\mathrm{const}.\notag
\end{align}
Assumptions (\ref{BETAPHI:1}) for the functions $r_{i1}$ are the
assumption that as $\varepsilon\to 0$ the functions
$r_{i}(t,\varepsilon)$ approximate continuous functions with a
possible discontinuity (jump) of the derivatives at $t=t^{*}$.

More precisely, the functions $r_i(t,\varepsilon)$ determined in
(\ref{PHI:1}) are approximations of the functions
\begin{equation*}
r_i(t,0)=r_{i0}+\psi_0 r^-_{i1} H(-\psi_0),
\end{equation*}
where $r^-_{i1}=\lim_{\tau\to-\infty}r_{i1}(\tau,t)$,
$r_i(t,\varepsilon)-r_i(t,0)=\mathcal{O}(\varepsilon)$, and $H(z)$
is the Heaviside function.

Now we describe the ansatz of the approximation of the temperature
$\bar{\sigma}(r,t)$. Here we also must use the continuation
procedure for matching the temperature for $t<t^{*}$, when there
are three subdomains of the domain~$Q$, and the temperature for
$t>t^{*}$, when there is only one phase. For this, we first
introduce a ``model'' of the temperature whose continuation is
reduced to the problem of continuation of functions depending only
on the time~$t$. This model must have the same structure as the
temperature $\bar{\sigma}$. Obviously, the simplest function of
this form is a function linear in~$x$ in the layers
$\Omega^{+}_{l,r}(t)$ and quadratic in the layer $\Omega^{-}(t)$.

Namely, we set
\begin{align*}
&\left.\frac{\partial \bar{\sigma}}{\partial r}
\right|_{r=\hat{r}_{1}\pm0}=\pm\gamma^{\pm}_{1},
\\
&\left.\frac{\partial \bar{\sigma}}{\partial r}
\right|_{r=\hat{r}_{2}\pm0}=\pm\gamma^{\pm}_{2},
\end{align*}

The functions $\gamma_{j}^{\pm}=\gamma_{j}^{\pm}(t)$, $j=1,2$,
depend on~$t$, are continuous for $t<t^{*}$, and have limits as
$t\to t^{*}-0$. Therefore, they can be continued to the interval
$[0,t_{1}]$, $t_{1}>t^{*}$, so that the properties of being
continuous and the signs are preserved. We shall use the previous
notation for the continued functions. Now we can introduce the
temperature model
\begin{align}\label{ModT:1}
\check{T}&=\gamma_{1}^{-}(t)(r_{1}-r)H(r_{1}-r)
+\gamma_{2}^{+}(t)(r-r_{2})H(r-r_{2})\notag\\
&+\gamma^{-}(r,t)\frac{(r_{1}-x)(x-r_{2})}{\psi}
H(r-r_{1})H(r_{2}-r)\\
&+\hat{\gamma}(r,t)\frac{(r_{1}-r)(r-r_{2})}{\psi}
H(r_{1}-r)H(r-r_{2})+I(r,t),\notag
\end{align}
where $\gamma^{-}(r,t)$, $\hat{\gamma}(r,t)$, and~$I(r,t)$ are
linear functions of~$r$,
$\left.\gamma^{-}\right|_{r=r_{i}}=\gamma^{-}_{i}(t)$,
$\left.\hat{\gamma}\right|_{r=r_{i}}=\hat{\gamma}_{i}(t)$, and
$\left.I\right|_{r=r_{i}}=(-1)^{i+1}(r_{i}r_{it}+\varkappa_{2})$.
In more detail, the functions $\hat{\gamma}(r,t)$ and
$\gamma^{-}(r,t)$ can be written in the form
\begin{align*}
&\gamma^{-}=\frac{\gamma_{1}^{+}-\gamma_{2}^{-}}{2}
-(r-r^{*})\frac{\gamma_{1}^{+}+\gamma_{2}^{-}}{\psi}
\\
&\hat{\gamma}=\frac{\hat{\gamma}_{1}+\hat{\gamma}_{2}}{2}
-(r-r^{*})\frac{\hat{\gamma}_{1}-\hat{\gamma}_{2}}{\psi}, \qquad
r^{*}=\frac{r_{1}+r_{2}}{2}.
\\
&I=\frac{k_{1}-k_{2}}{2} -(r-r^*)\frac{k_{1}+k_{2}}{\psi}, \quad
k_{1}=r_{1}r_{1t}+\varkappa_{2},\quad
k_{2}=r_{2}r_{2t}+\varkappa_{2}.
\end{align*}
We have the relations
\begin{equation}\label{OPhiPhi}
r_{i}(t,\varepsilon)=\hat{r}_{i}(t) +\mathcal{O}(\varepsilon),
\qquad i=1,2.
\end{equation}
This implies that, for $t<t^{*}$,
\begin{align}\label{TDelta}
&\left.\left[\frac{\partial \check{T}}{\partial r}\right]
\right|_{r=r_{i}}= \left.\left[\frac{\partial
\bar{\sigma}}{\partial r}\right]
\right|_{r=\hat{r}_{i}},\\
&\left[\frac{\partial \check{T}}{\partial r}\right]\delta(r-r_{i})
=\left[\frac{\partial \bar{\sigma}}{\partial r}\right]
\delta(r-\hat{r}_{i})+\mathcal{O}_{\mathcal{D}'}(\varepsilon).\notag
\end{align}
We rewrite the product
\begin{equation*}
H(r-r_{1})H(r_{2}-r)
\end{equation*}
as follows
\begin{gather}\label{ProdH:1}
H(r-r_{1})H(r_{2}-r) =B(\rho)[H(r-r_{1})-H(r-r_{2})]
+\mathcal{O}_{\mathcal{D}'}(\varepsilon).
\\
H(r_{1}-r)H(r-r_{2}) =(1-B(\rho))[H(r-r_{2})-H(r-r_{1})]
+\mathcal{O}_{\mathcal{D}'}(\varepsilon). \notag
\end{gather}
Here $B(\rho)\in C^\infty$, $B(\rho)\to1$, $\rho\to\infty$,
$B(\rho)\to0$, $\rho\to-\infty$,
$B^{(\alpha)}=\mathcal{O}(|\rho|^{-N})$, $|\rho|\to\infty$ for any
$N>0$, $\alpha>0$.

In what follows, we write $B(\tau)$ instead of $B(\rho)$. The
difference is that $B(\rho)$ is an unknonwn function ($\rho$ is
unknown), while $B(\tau)$ is a known function.

But we do not use this remark and, by definition, we set
\begin{align*}
\check{T}&=I+\left[-\gamma_{1}^{+}(r_{1}-r)+
\frac{\gamma^{-}}{\psi}(r_{1}-r)(r-r_{2})B(\rho)-\right.\\
&\left.\frac{\hat{\gamma}}{\psi}(r_{1}-r)(r-r_{2})\left(1-B(\rho)\right)\right]H(r-r_{1})+\\
&\left[\gamma_{2}^{+}(r_{1}-r)-
\frac{\gamma^{-}}{\psi}(r_{1}-r)(r-r_{2})B(\rho)+\right.\\
&\left.\frac{\hat{\gamma}}{\psi}(r_{1}-r)(r-r_{2})\left(1-B(\rho)\right)\right]H(r-r_{2}).
\notag
\end{align*}

Now we can write the expression $\partial^{2}\check{T}/\partial
r^{2}$ uniformly in $\psi=r_{2}-r_{1}$. We have
\begin{align}\label{D2T}
\frac{\partial^{2}\check{T}}{\partial x^{2}}&=
\frac{2}{\psi}\left[-\gamma^{-}_{1}(t)B(\rho)+\hat{\gamma}_{1}(t)\left(1-B(\rho)\right)\right]H(r-r_{1})+\notag\\
&\frac{2}{\psi}\left[\gamma^{-}_{2}(t)B(\rho)+\hat{\gamma}_{2}(t)\left(1-B(\rho)\right)\right]H(r-r_{2})+\\
&\left[\gamma^{+}_{1}(t)+\gamma^{-}_{1}(t)B(\rho)-\hat{\gamma}_{1}(t)\left(1-B(\rho)\right)\right]\delta(r-r_{1})+\notag\\
&\left[\gamma^{+}_{2}(t)+\gamma^{-}_{2}(t)B(\rho)-\hat{\gamma}_{2}(t)\left(1-B(\rho)\right)\right]\delta(r-r_{2})\notag,
\end{align}
We set $\gamma_{1}^{+}=-\hat{\gamma}_{1}$ and
$\gamma_{2}^{-}=\hat{\gamma}_{2}$ and see that the coefficients
for $\delta(x-r_{i})$ in (\ref{D2T}) take the form
\begin{align}\label{CoefDelta}
&\left[\gamma_{1}^{-}+\gamma_{1}^{+}\right]B(\tau)\delta(r-r_{1})
+\left[\gamma_{2}^{-}+\gamma_{2}^{+}\right]B(\tau)\delta(r-r_{1}).
\end{align}
We note that the replacement of $B(\rho)$ by $B(\tau)$ implies a
correction of order $\mathcal{O}_{\mathcal{D}'}(\varepsilon)$. Now
we note that the inequalities
\begin{align*}\label{HHPsi}
\left[H(r-r_{1})-H(r-r_{2})\right]/\psi >0,\qquad
\left[H(r-r_{2})-H(r-r_{1})\right]/\psi <0,
\end{align*}
hold for $\psi\neq 0$. Also, we note that
$\partial\check{T}/\partial r$ is a smooth function. So, by
(\ref{TDelta}) and the Stefan condition (\ref{SSt:1}) we have
\begin{equation*}
\gamma_{1}^{-}+\gamma_{1}^{+}=r_{10}^{3}(t)r_{10}'(t),\qquad
\gamma_{2}^{-}+\gamma_{2}^{+}=-r_{20}^{3}(t)r_{20}'(t)
\end{equation*}
for $t<t^*$ and assume that the continuations of the functions
contained in these relations are chosen so that these relations
hold for $t\in [0,t_1]$.

The assumption that $\overline{\sigma}(r,t)$ is the classical
solution of problem (\ref{SH:1})--(\ref{SSt:1}) implies that
\begin{equation}\label{DifferTS}
\bar{\sigma}-\check{T} =\check{q}\in C([0,t^{*});C^{1}(\Omega))
\end{equation}
and the limit $\check{q}$ exists as $t\to t^{*}-0$ and for
$x\in\Omega$. Therefore, we can continue the function
$\check{q}$ to the domain $\Omega\times[0,t_{1}]$, where
$t_{1}>t^{*}$ and the properties of smoothness are preserved,
and moreover,
\begin{equation*}
\left.\check{q}\right|_{r=\hat{r}_{i}}=0, \qquad i=1,2, \qquad
t<t^{*}.
\end{equation*}
We shall construct the global temperature $\check{\theta}$
in the form
\begin{equation}\label{CheckS}
\check{\sigma}=e(r)\check{T}+\check{q}+\hat{q},
\end{equation}
where $\hat{q}$ is the desired function, $e(r)\in
C_{0}^{\infty}([R_{1},R_{2}])$, and $e\equiv 1$ for
$r\in[\hat{r}_{1}(0),\hat{r}_{2}(0)]$.

\section{Construction of the weak asymptotic
solution}\label{Constr:1}

We first consider the heat equation (\ref{PhHeat}). Using the
formulas of the weak asymptotics method and taking (\ref{AU:1})
into account, we obtain
\begin{align}\label{Udivt:1}
\frac{\partial \check{u}}{\partial t}&=2\pi
r_{1}^{2}\left[r_{1t}(2-B_{\dot{0}0})+
\frac{\beta_{\tau}\psi_{0}'}{\beta^{2}}
B^{z}_{\dot{0}0}\right]\delta(r-r_{1})\\
&+2\pi r_{2}^{2}\left[-r_{2t}(2-B_{\dot{0}0})+
\frac{\beta_{\tau}\psi_{0}'}{\beta^{2}}
B^{z}_{\dot{0}0}\right]\delta(r-r_{2})+
\mathcal{O}_{\mathcal{D}'}(\varepsilon),\notag
\end{align}
where the estimate $\mathcal{O}_{\mathcal{D}'}(\varepsilon)$ is
uniform in $\tau$, $\psi_{0}=r_{20}-r_{10}$,
\begin{equation*}
B_{\dot{0}0}=\int_{\mathbb{R}^{1}}\dot{\omega}_{0}(z)
\omega_{0}(-\eta-z)\,dz, \quad
B^{z}_{\dot{0}0}=\int_{\mathbb{R}^{1}}
z\dot{\omega}_{0}(z)\omega_{0}(-\eta-z)\,dz,\quad \eta=\rho\beta.
\end{equation*}

We assume that function $\check{q}$ and its continuation are
chosen so that the function $\check{q}$ satisfies the boundary
conditions of the original problem.
Then the boundary conditions for the function $\hat{q}$ are
zero.

We substitute the function $\check{\sigma}$ determined by relation
(\ref{CheckS}) and expression (\ref{Udivt:1}) for
$\partial\check{u}/\partial t$ into Eq.~(\ref{PhHeat}). We
denote
\begin{equation*}
\hat{L}=\frac{\partial}{\partial t}-\frac{\partial^{2}}{\partial
r^{2}}.
\end{equation*}
According to the definition~\ref{Def:1} with accuracy
$\mathcal{O}_{\mathcal{D}'}(\varepsilon)$, we obtain
\begin{equation*}
\hat{L}\check{\sigma}+r\frac{\partial \check{u}}{\partial t}
=\hat{L}e\check{T}+\hat{L}q+A_1\delta(r-r_1) +A_2\delta(r-r_2),
\end{equation*}
where
\begin{equation*}
q=\check{q}+\hat{q},\qquad A_i=
\frac{r_{i}^{3}}{2}\left[(-1)^{i+1}r_{it}(2-B_{\dot{0}0})
+\frac{\beta_\tau\psi'_{0}}{\beta^2}B^z_{\dot{0}0}\right],\qquad
i=1,2.
\end{equation*}
We let $\tau$ tend to $\infty$ (i.e., for $t<t^*$) and, in view of
(10) and (18), obtain
\begin{equation*}
\hat{L}(e\check{T}+{q})=\hat{r}_{2}^{3}r_{2}'
\delta(r-\hat{r}_{2}) - \hat{r}_{2}^{3}r_{2}'
\delta(r-\hat{r}_{1}).
\end{equation*}

Since $e\check{T}+\check{q}\to\overline{\sigma}$ as
$\tau\to\infty$, we see that $\hat{q}\to0$ as $\tau\to\infty$. The
total equation for the function $q=\check{q}+\hat{q}$ has the form
\begin{align}\label{Lq}
&\hat{L}(q+e(x)I) = F(x,t,\tau)\notag
\\
&\qquad -\frac{\partial^2}{\partial r^2}
\bigg(\gamma^-\frac{(r_1-r)(r-r_2)}{\psi}\bigg)
B(\tau)[H(r-r_{1})-H(r-r_{2})]\\
&\qquad -\frac{\partial^2}{\partial r^2}
\bigg(\hat\gamma\frac{(r_1-r)(r-r_2)}{\psi}\bigg) (1-B(\tau))
[H(r-r_{2})-H(r-r_{1})], \quad r\in\Omega, \notag
\end{align}
where $F(r,t,\tau)$ is a piecewise continuous function containing
of terms that are uniformly bounded in $\varepsilon$ on
$\Omega\times[0,t_1]$.

Equation (\ref{Lq}) was derived under the assumption that
\begin{align}\label{SumD}
\sum^{2}_{i=1}\{B[\gamma^+_i+\gamma^-_i]-A_i\}\delta(r-r_i)
=\mathcal{O}_{\mathcal{D}'}(\varepsilon).
\end{align}
By $\hat{q}_i$, $\hat{q}^*_i$
we denote the terms in $q$ such that
\begin{align}\label{qqqq}
\hat{L}\hat{q}_1&=-\frac{2(\gamma^+_1+\gamma^-_2)}{\psi}B(\tau)
[H(r-r_2)-H(r-r_1)],
\\
L\hat{q}_2&=\frac{2(\gamma^-_1+\gamma^+_2)}{\psi}(1-B)
[H(r-r_2)-H(r-r_1)],
\notag\\
\hat{L}\hat{q}^*_1&=2(r-r^*)(\gamma^+_1-\gamma^-_2)\psi^{-2}
B(\tau) [H(r-r_2)-H(r-r_1)],
\notag\\
\hat{L}\hat{q}^*_2&=-2(r-r^*)(\gamma^-_1-\gamma^+_2)\psi^{-2}
B(\tau) [H(r-r_2)-H(r-r_1)].\notag
\end{align}

With accuracy up to functions smooth in $\Omega$,
we can calculate them using the fundamental solution
of the heat equation (in what follows, we shall consider
only the first equation),
\begin{align}\label{hatq}
\hat{q}_1&=-\frac{1}{2\sqrt{2\pi}}\int^{t}_{0}
\frac{(\gamma^+_1(\alpha)+\gamma^-_2(\alpha))}{\sqrt{t-\alpha}}
B(\tau(\alpha))\\
&\times \int^{r_2}_{r_1} \frac{ \exp\{-(r-\xi)^2/(t-\alpha)\}
}{\psi}\,d\alpha d\xi.\notag
\end{align}
It is clear that the functions $\hat{q}_i$ are uniformly bounded.
It is also clear that $\hat{q}_i\not\in C^1(\Omega)$, but
$\hat{q}_i\in C^{1,2}(\Omega\times[0,t_1])\setminus
(\{r=r_1\}\cup\{r-r_2\})$. Therefore, to justify Eq.~(\ref{Lq}),
we must verify that no $\delta$-functions arise in calculating the
derivatives $\partial^2 \hat{q}_i/\partial r^2$, $\partial^2
\hat{q}^*_i/\partial r^2$ and $\partial \hat{q}_i/\partial t$,
$\partial \hat{q}^*_i/ \partial t$.

For this, we note that, for any test function $\zeta(r)$ up to
functions smooth in $\Omega$, we have the relation
\begin{equation}\label{31}
\langle \hat{q}_i,\zeta(r)\rangle
=\frac12\int^{t}_{0}\gamma^-_1(\alpha) B\int^{r_2}_{r_1}\psi^{-1}
f(\xi,t-\alpha)\,d\alpha]\,d\xi,
\end{equation}
where $f(\xi,t-\alpha) =\int_{\mathbb{R}^{1}}
\left[\zeta(r)e^{-(r-\xi)^2/(t-\alpha)}\right]/\left[\sqrt{2\pi(t-\alpha)}\right]dr$
is the solution of the heat equation $f_t-f_{\xi\xi}=\zeta(\xi)$
at the point $t-\alpha$. It is clear that $f(\xi,t)\in
C^\infty(\Omega)$ for all~$t$. Calculating the derivative
$\partial \hat{q}_1/\partial r$, we obtain
\begin{equation*}
\left\langle \frac{\partial\hat{q}_1}{\partial
r},\zeta\right\rangle =\langle \hat{q}_1,\zeta'\rangle
=2\int^{t}_{0}\gamma^-_1(\alpha)B \int^{r_2}_{r_1}\psi^{-1}
f_1(\xi,t-\alpha)\,d\alpha\,d\xi,
\end{equation*}
where $f_1(\xi,t)$ is a solution of the equation
$f_t-f_{\xi\xi}=\zeta''(\xi)$. It is clear that the relation
$$
\left\langle \frac{\partial^2\hat{q}_i}{\partial
r^2},\zeta(r)\right\rangle \not= \sum
G_i\zeta(r_i)+\mathcal{O}(\varepsilon).
$$
cannot hold for any coefficients $G_i$. Similarly to (\ref{31}),
with accuracy up to smooth functions, we can write the functions
$\hat{q}^*_i$, $i=1,2$. They have the same properties as
$\hat{q}_i$ with the only additional condition
$$
\hat{q}^0_i|_{r=r^*}=0.
$$
This easily follows from the explicit formula of the type of (31)
and the fact that, for $r=r^*$, the integral with respect to~$\xi$
is an integral of an odd function over a symmetric interval.
Therefore, condition (\ref{hatq}) is necessary for deriving
Eq.~(\ref{Lq}). This implies that the function $\hat{q}+eI$ is
uniformly bounded in $\Omega\times[0,t_1]$ and belongs
to~$C^{1,0}$ for $r\not=r_i$, $i=1,2$.

So if the function $\hat{q}+eI$ satisfies Eq.~(\ref{Lq}) with zero
initial conditions and the boundary conditions that follow from
the fact that the function
$\check{\theta}=e\check{T}+\check{q}+\hat{q}$ must satisfy the
boundary conditions of the original problem, while the functions
$\check{T}$ and $\check{q}$ are known (more precisely, they will
be determined after the functions $\varphi_i$ are found), then we
have the following assertion.

\begin{theorem}\label{TH:1}
Suppose that the function $\check{\sigma}$ is determined by
relation (\ref{CheckS}), and the function $\hat{q}$ is a solution
of Eq.~(\ref{Lq}) with zero initial condition and zero boundary
conditions. Suppose that relation (\ref{SumD}) holds.

Then the pair of functions $\check{\sigma}$ and $\check{u}$ is a
weak asymptotic solution of the heat equation in the phase field
system, i.e., relation (\ref{PhHeat}) holds.
\end{theorem}

\begin{remark}\label{Rem:1}
Expression (\ref{SumD}) may seem to be strange, because we did not
take into account the well-known fact that the Dirac functions are
linearly independent. But this fact is taken into account in the
framework of the weak asymptotics method (see Lemma~2, Section~5).
Here we do not want to mix the substitution of the ansatz into the
equation and the analysis of the results of this substitution (see
Section~4).
\end{remark}

Now we consider the Allen--Cahn equation (\ref{PhAC}). We must
calculate the weak asymptotics of the expression
\begin{equation*}
\mathcal{F}\buildrel\rm def\over=\frac{\partial
\check{u}}{\partial r}\left[\varepsilon
L\check{u}-\frac{\left(\check{u}-\check{u}^{3}\right)}{\varepsilon}
-\frac{\check{\sigma}}{r}\right],
\end{equation*}
where the function $\check{\sigma}$ is determined by the relation
(\ref{CheckS}). We obtain
\begin{align}\label{FCal:1}
\mathcal{F}=&V_{1}^{1}\delta(x-\varphi_{1})+V_{2}^{1}\delta(x-\varphi_{2})
\\
&+V_{1}^{2}\delta'(x-\varphi_{1})+V_{2}^{2}\delta'(x-\varphi_{2})
+\mathcal{O}_{\mathcal{D}'}(\varepsilon),\notag
\end{align}
where $V_{j}^{i}$, $i,j=1,2$, are linear combinations of several
convolutions. Their expressions will be given below. Here we note
the following. We have
\begin{align}
\varepsilon
&\int_{\mathbb{R}^{1}}\check{u}_{r}(r^{2}\check{u}_{r})_{r}\zeta
dr+\frac{1}{\varepsilon}\int_{\mathbb{R}^{1}}\check{u}_{r}\left(\check{u}-\check{u}^{3}\right)\zeta
dr=\label{ururr}\\
&\frac{\varepsilon}{2}\int_{\mathbb{R}^{1}}r\check{u}^{2}_{r}\zeta
-\int_{\mathbb{R}^{1}}r^{2}\left(\frac{1}{\varepsilon}F(\check{u})+
\frac{\varepsilon}{2}\check{u}^{2}_{r}\right)\zeta_{r}dr,\notag
\end{align}
where
\begin{equation*}
F(\check{u})=\frac{\check{u}^{4}}{4}-\frac{\check{u}^{2}}{2}+\frac{1}{4}.
\end{equation*}

We denote
\begin{equation}\label{Omega}
\Omega(z,\eta)=\frac{1}{2}\left\{1+\omega_{0}(z)+\omega_{0}(-z-\eta)
-\omega_{0}(z)\omega_{0}(-z-\eta)\right\}.
\end{equation}

We have the following estimates:
\begin{align}
\Omega(z,\eta)&=1+f(z,\eta) e^{2\eta}, \qquad
\eta\to-\infty,\label{r34}
\\
\Omega(z,\eta)&=\omega_0(z)+f_1(z,\eta) e^{-2\eta}, \qquad
\eta\to\infty,\label{r35}
\end{align}
where
$$
\int|f(z,\eta)|\,dz\leq\mathrm{const},\qquad
\int|f_1(z,\eta)|\,dz\leq\mathrm{const}.
$$

These relations readily follow from the explicit form
of the function $\omega_0(z)$,
$$
\omega_0(z)=\frac{e^z-e^{-z}}{e^z+e^{-z}}.
$$
In fact, relations (\ref{r34}) and (\ref{r35}) express the
above-described properties of the ansatz $\check{u}$~(\ref{AU:1})
in different terms.

Using the technique of the weak asymptotics method (see
Lemma~\ref{lemma7} in Section~5), we can write
\begin{align}\label{urut}
\int_{\mathbb{R}^1}r^{2}\check u_r \check u_t\zeta dr
&=\sum^2_{i=1}r^{2}_{i}r_{it}
\int\Omega'_z(\Omega'_z-\Omega'_\eta)\,dz
\\
&\qquad +r^{2}_{i} \frac{\beta_\tau\psi'_{0}}{\beta^2}
\int\Omega'_z(\Omega'_z-\Omega'_\eta)\,dz
+\mathcal{O}_{\mathcal{D}'}(\varepsilon). \notag
\end{align}
Similarly, we obtain (see Lemmas~\ref{lemma5} and~\ref{lemma7},
Section~5):
\begin{align}\label{ThUX}
&\int_{\mathbb{R}^{1}}r\check{u}_{r}\check{\sigma}\zeta dr
=-r_{1}\left.\left(eI+\hat{q}+\check{q}\right) \right|_{r=r_{1}}
\int_{\mathbb{R}^{1}}
\dot{\Omega}_{\eta}(z,\eta)\zeta(r_{1})dz\\
&+r_{2}\left.\left(eI+\hat{q}+\check{q}\right) \right|_{r=r_{2}}
\int_{\mathbb{R}^{1}}\dot{\Omega}_{\eta}(z,\eta) \zeta(r_{2})\,dz
+\mathcal{O}(\varepsilon^\mu),\quad \mu\in(0,1/2).\notag
\end{align}

Thus, adding expressions (\ref{ururr}), (\ref{urut}) and
(\ref{ThUX}), we see that the coefficients of the
$\delta$-functions in formula (\ref{FCal:1}) have the form
\begin{equation}\label{VVV}
V^1_i=-r_{i}^{2}\left[r_{it}B_\Omega
+\frac{\beta_\tau\psi'_{0}}{\beta^2} B^z_\Omega\right]
+(-1)^{i+1}r_{i}\left.\left(eI+q\right)\right|_{r=r_i}C_\Omega-r_{i}\hat{C},
\end{equation}
where
\begin{align}\label{BVVV}
i=1,2,\quad B_\Omega&=\int\Omega'_z(\Omega'_z-\Omega'_\eta)\,dz,
\notag\\
B^z_\Omega&=\int[z(\Omega'_z-\Omega'_\eta)-(z+\eta)\Omega'_\eta]
(\Omega'_z-\Omega'_\eta)\,dz,
\\
C_\Omega&=\int(\Omega'_z-\Omega'_\eta)\,dz,\qquad
q=\hat{q}+\check{q} \notag
\end{align}
and
\begin{equation}\label{C:1}
\hat{C}=\frac{1}{2}\int_{\mathbb{R}^{1}}(\Omega'_{z})^{2}\,dz
\end{equation}

Similarly, we obtain (see Lemma~\ref{lemma6}, Section~5):
\begin{equation}\label{V12:1}
V_{1}^{2}=r_{1}^{2}\left(\beta\hat{C}+\frac{1}{\beta}\hat{D}\right),
\quad
V_{2}^{2}=r_{2}^{2}\left(\beta\hat{C}+\frac{1}{\beta}\hat{D}\right),
\end{equation}
where
\begin{equation}\label{D:1}
\hat{D}=-\frac{1}{2}\int_{\mathbb{R}^{1}} F(\Omega)\,dz
\end{equation}

Thus, to obtain $\beta=\beta(\eta)$
we have the equation
(obviously, this is a necessary condition for the relation
$\mathcal{F}=\mathcal{O}_{\mathcal{D}'}(\varepsilon)$ to hold,
see~(\ref{FCal:1})):
\begin{equation}\label{BETA:1}
\beta^{2}=\frac{\hat{D}}{\hat{C}}
\end{equation}
According to (\ref{C:1}) and (\ref{D:1}) ($\hat{D}$ is positive),
the right-hand side of the relation (\ref{BETA:1}) is positive.

So we have proved the following assertion.

\begin{theorem}\label{TH:2}
Suppose that the assumptions of Theorem~{\rm1} and (\ref{BETA:1})
are satisfied and
\begin{equation}\label{SumVDelta:1}
\sum^{2}_{i=1}V^1_1\delta(x-\varphi_1)
+V^1_2\delta(x-\varphi_2)=O_{D'}(\varepsilon).
\end{equation}
Then the pair of functions $\check\sigma$, $\check u$ is a weak
asymptotic solution of the phase field system (\ref{PhHeat}),
(\ref{PhAC}) in the sense of Definition~\ref{Def:1}.
\end{theorem}

Thus, under the assumption that the classical solution of the
phase field system exists (see~(i) and~(ii) above), relations
(\ref{Lq}), (\ref{SumD}), (\ref{BETA:1}), and (\ref{SumVDelta:1})
are sufficient conditions for constructing a weak asymptotic
solution of system (\ref{PhHeat}), (\ref{PhAC}). In what follows,
we prove that the relations mentioned above are equations for
determining the functions $\beta=\beta(\tau)$ and $r_{i1}(\tau)$,
$i=1,2$. We present an algorithm for solving these equations.

\section{An analysis of the relations obtained}

We begin with relation (\ref{SumD}) and, for a while, forget
everything said about the notion of linear independence.

Then, for this relation to hold, it suffices to have the two
relations
\begin{equation}\label{BA}
B[\gamma^+_i+\gamma^-_i]=A_i,\qquad i=1,2.
\end{equation}
In view of (\ref{TDelta}), we have $\gamma^+_i+\gamma^-_i
=\left.\left[\frac{\partial\overline{\sigma}}{\partial r}\right]
\right|_{r=\hat{r}_i(t)}$. Taking this into account and adding
relations (\ref{BA}) and Lemma~\ref{Lemma:2}, we obtain
\begin{align}\label{I1}
&\psi'_0\left(\rho_\tau-\frac12\int_{\mathbb{R}^{1}}
\dot\omega(z)\omega_0(-z-\eta)
\left(\rho_\tau-2\frac{\beta_\tau}{\beta}z\right)\right)\,dz\\
&\qquad =B(\tau)\left\{
\left.\left[\frac{\partial\overline{\theta}}{\partial r}\right]
\right|_{r=\hat{r}_1(t)}
+\left.\left[\frac{\partial\overline{\theta}}{\partial r}\right]
\right|_{r=\hat{r}_2(t)} \right\} =2B(\tau)\psi'_{0}.\notag
\end{align}
Here we used the relation
$$
\psi_t=r_{2t}-r_{1t} =\psi'_0\rho_\tau+\psi_0(r_{21}-r_{11})_t,
$$
which follows from the definition of the function $\rho$,
see~(\ref{R1PH1}).

We set the last term to be zero, since we show below that
$r_{21}-r_{11} =\mathcal{O}(|\tau|^{-1})$. In view of Lemma~3 in
Section~5, in this case we have
$(r_{21}-r_{11})_t\psi_0=\mathcal{O}(\varepsilon)$. Moreover, in
view of the Stefan conditions and the choice of the continuation
of the functions $\gamma^\pm_i$ and $r_{0i}(t)$, we have
$$
\gamma^+_1+\gamma^+_2+\gamma^-_1+\gamma^-_2
=\left.\left[\frac{\partial\overline{\theta}}{\partial r}\right]
\right|_{r=\hat{r}_1(t)}
+\left.\left[\frac{\partial\overline{\theta}}{\partial r}\right]
\right|_{r=\hat{r}_2(t)} =2\psi'_{0}.
$$

We transform the left-hand side of relation (\ref{I1}):
\begin{align*}
I&\buildrel{\rm def}\over=\psi_{0}'\left[\rho_{\tau}
-\frac{1}{2}\int_{\mathbb{R}^{1}}
\dot{\omega}_{0}(z)\omega_{0}(-z-\eta)
\left(\rho_{\tau}-\frac{2\beta_{\tau}}{\beta^{2}}z\right)dz\right]\\
&=\psi_{0}'\left[\rho_{\tau}+ \frac{1}{2}\int_{\mathbb{R}^{1}}
\dot{\omega}_{0}(z)\omega_{0}(-z-\eta)
\frac{\beta_{\tau}}{\beta^{2}}z\,dz\right.\\
&\left. -\frac{1}{2}\int_{\mathbb{R}^{1}}
\dot{\omega}_{0}(z)\omega_{0}(-z-\eta)
\left(\rho_{\tau}+\frac{\beta_{\tau}}{\beta^{2}}(-z-\eta)+
\frac{\beta_{\tau}}{\beta^{2}}\eta\right)dz\right].
\end{align*}
Next, we change the variables in the last integral
\begin{equation*}
\frac{\beta_{\tau}}{\beta^{2}}
\int_{\mathbb{R}^{1}}(-z-\eta)\dot{\omega}_{0}(z)
\omega_{0}(-z-\eta)\,dz
=-\frac{\beta_{\tau}}{\beta^{2}}
\int_{\mathbb{R}^{1}}
z\dot{\omega}_{0}(-z-\eta)\omega_{0}(z)\,dz
\end{equation*}
and note that
\begin{equation*}
\rho_{\tau}-\frac{\beta_{\tau}}{\beta^{2}}\eta
=\frac{1}{\beta}\frac{\partial}{\partial \tau}\eta.
\end{equation*}
Finally, we get
\begin{align*}
I&=\psi_{0}'\left[\rho_{\tau}- \frac{1}{2}\int_{\mathbb{R}^{1}}
\left\{z\dot{\omega}_{0}(z)\omega_{0}(-z-\eta)
+z\dot{\omega}_{0}(-z-\eta)\omega_{0}(z)\right\}\,dz\right.\\
&\qquad -\left.\frac{1}{2\beta}\frac{\partial}{\partial
\eta}\int_{\mathbb{R}^{1}}
\dot{\omega}_{0}(z)\omega_{0}(-z-\eta)\,dz\right].
\end{align*}
Since the function $\omega_{0}(z)$ is odd, we see that
the expression in braces
in the first integral in the right-hand side  is
\begin{equation*}
-z\frac{\partial}{\partial z}
\left(1-\omega_{0}(z)\omega_{0}(z+\eta)\right).
\end{equation*}
Hence, integrating by parts, we obtain
\begin{align*}
I&=\psi_{0}'\left[\rho_{\tau}-\frac{\beta_{\tau}}{2\beta^{2}}
\int_{\mathbb{R}^{1}}
\left(1-\omega_{0}(z)\omega_{0}(z+\eta)\right)\,dz\right.\\
&\qquad -\left.\frac{1}{2\beta}\frac{\partial\eta}{\partial\tau}
\int_{\mathbb{R}^{1}}
\dot\omega_{0}(z+\eta)\omega_{0}(z)\,dz\right].
\end{align*}
Or, finally,
\begin{equation}\label{I:1}
I=\psi_{0}'\left[\rho_{\tau}
+\frac{1}{2}\frac{\partial}{\partial\tau}\left(\frac{1}{\beta}\tilde{B}(\eta)\right)\right],
\end{equation}
where
\begin{equation*}
\tilde{B}(\eta)
=\int_{\mathbb{R}^{1}}
\left(1-\omega_{0}(z+\eta)\omega_{0}(z)\right)\,dz
=2\eta\tanh\eta.
\end{equation*}
So, in view of (\ref{I1}), we obtain the equation for $\rho$:
\begin{equation}\label{RhoSt:1}
\frac{\partial}{\partial\tau}\left(\rho
+\frac{1}{2\beta}\tilde{B}(\eta)\right)=2B(\tau).
\end{equation}
Or
\begin{equation*}
\frac{\partial}{\partial\tau}
\left(\frac{\eta+\eta\tanh\eta}{\beta}\right)=2B(\tau).
\end{equation*}
After integration, we obtain
\begin{equation}\label{ETA1}
\eta(1+\tanh\eta)= 2\beta\int^{\tau}_{-\infty}B(\alpha)\,d\alpha.
\end{equation}

Since $B(\tau)\to1$ as $\tau\to\infty$,
we have $\eta\tau^{-1}\to1$ as $\tau\to\infty$.

Since $B(\tau)=\mathcal{O}(|\tau|^{-N})$ for any~$N$ as
$\tau\to-\infty$, the integral in the right-hand of (\ref{ETA1})
converges and hence the left-hand side tends to the limit
\begin{equation*}
\lim_{\tau\to-\infty}\eta(1+\tanh\eta)
=\lim_{\tau\to-\infty}2\beta\int^{\tau}_{-\infty}B(\alpha)\,d\alpha=0.
\end{equation*}
Moreover, in view of the inequality $B\geqslant0$, we have
$\eta\geqslant0$ for $\tau\in\mathbb{R}^2_\tau$, which implies
that
\begin{equation*}
\eta\to0\qquad\text{as}\quad \tau\to-\infty.
\end{equation*}
Here we took into account the inequality $\beta>0$, which follows
from~(\ref{BETA:1}).

Substituting $\beta=(\hat{C}\hat{D}^{-1}(\eta))^{1/2}$ into
(\ref{ETA1}), we obtain the following equation for the function
$\eta$:
\begin{equation}\label{ETA2}
\eta(1+\tanh\eta)
=(\hat{C}\hat{D}^{-1})^{1/2}\int^{\tau}_{0}B\,d\tau'.
\end{equation}

Since the functions contained in (\ref{ETA2}) are monotone and the
limits exist as $\tau\to\pm\infty$, this equation is obviously
solvable.

Next, we have
\begin{equation*}
\beta(\tau)=\sqrt{\frac{\hat{C}(\eta)}{\hat{D}(\eta)}}\to\beta^-
=\mathrm{const},\qquad \tau\to-\infty.
\end{equation*}
Moreover, as readily follows from the exponential rate of
convergence to the limit of the functions $\hat{C}/\hat{D}$,
$1-B$, and $1+\tanh\eta$ as $\tau\to\pm\infty$, the derivatives
satisfy the estimate
\begin{equation*}
\frac{\partial^\alpha\beta}{\partial\tau^\alpha}
=\mathcal{O}(|\tau|^{-N}),\qquad |\tau|\to\infty.
\end{equation*}
Now we calculate the limits of the expressions
\begin{equation*}
J_i=B(\gamma^+_i+\gamma^-_i)-A_i,\qquad i=1,2,
\end{equation*}
as $\tau\to\pm\infty$. As $\tau\to\infty$, we have
\begin{equation*}
\rho\to\infty\quad(\rho\sim\tau),\qquad B\to1,\qquad
A_i\to(-1)^{i+1}r_{i0}^3r_{i0}'.
\end{equation*}
Therefore, we have
\begin{equation*}
\lim_{\tau\to\infty}J_i=0
\end{equation*}
in view of the Stefan conditions (\ref{SSt:1}), and
$J_i=\mathcal{O}(\tau^{-1})$, $i=1,2$.

As $\tau\to-\infty$, we proved that $\eta\to0$
and hence
\begin{equation}\label{52}
r_1-r_2=\mathcal{O}(\varepsilon),\qquad \tau\to-\infty.
\end{equation}

\begin{theorem}
Relation (\ref{ETA2}) is a sufficient condition for relation
(\ref{SumD}) to hold.
\end{theorem}
\begin{proof} In view of the corollary of Lemma~\ref{Lemma:2}, Section~5,
the estimate for the coefficients of the $\delta$-functions
in~(\ref{SumD}) as $\tau\to\infty$ and (\ref{52}) are sufficient
conditions for relation (\ref{SumD}) to follow from (\ref{I1}). In
turn, relation (\ref{ETA2}) follows from (\ref{I1}).
\end{proof}

Now we analyze the relations $V^2_i=0$ and (\ref{SumVDelta:1}). We
use the explicit form of the function $F(z)$ and the formula
\begin{equation}
\omega_0=\frac{e^z-e^{-z}}{e^z+e^{-z}}.
\end{equation}

As $\eta\to\infty$, we have
\begin{equation}\label{54}
\Omega(z,\eta)=\omega_0(z)+\mathrm{O}(e^{-2\eta})
\end{equation}
and hence $\beta\to1$ as $\eta\to\infty$.

From relation (\ref{Omega}) we easily obtain
\begin{gather*}
\Omega'_z(z,\eta)=\mathcal{O}(e^{-|z|}),\qquad |z|\to\infty,\\
\Omega'_\eta(z,\eta)=\mathcal{O}(e^{-|z|}),\qquad |z|\to\infty,\\
C_\Omega=\frac{1}{4}\int_{\mathbb{R}^{1}} \dot\omega_0(z)(1-\omega_0(-z-\eta))\,dz\geq0.
\end{gather*}

Similarly to (\ref{54}), using (\ref{BA}), we can verify the
relations
\begin{align}\label{BC}
&\beta_\tau=\mathcal{O}(e^{-2|\eta|}), \qquad  \eta\to\infty,
\notag\\
&B_\Omega= 1+\mathcal{O}(e^{-2|\eta|}),\quad\quad\eta\to\infty,
\\
&C_\Omega=4+\mathcal{O}(e^{-2\eta}),\qquad \eta\to\infty,
\notag\\
&B^z_\Omega=0,\qquad \eta\to\infty. \notag
\end{align}

We rewrite the expressions for the coefficients $V^1_i$ in more detail:
\begin{align}
V^1_1&=-r_{1}^{2}r_{1t}B_\Omega
-\frac{r_{1}^{2}\beta'_\tau\psi'_{0}}{\beta^2}B^z_\Omega+r_{1}\left.\left(eI+q\right)\right|_{r=r_1}
C_\Omega-r_{1}\hat{C},\label{V11}
\\
V^1_2&=-r_{2}^{2}r_{2t}B_\Omega
-\frac{r_{2}^{2}\beta'_\tau\psi'_{0}}{\beta^2}B^z_\Omega-r_{2}\left.\left(eI+q\right)\right|_{r=r_2}
C_\Omega-r_{2}\hat{C}.\label{V12}
\end{align}
We take the limit in formulas (\ref{V11}), (\ref{V12}) as $\tau\to
\infty$ i.e. prior to interaction. So, in view of conditions
(\ref{C:1}), (\ref{BC}) we obtain
\begin{align*}
&\left.\frac{\bar{\sigma}}{r}\right|_{r=r_{10}(t)}=\hat{C}^{+}r_{10}'(t)+\frac{\hat{C}^{+}}{r_{10}(t)},\\
&\left.\frac{\bar{\sigma}}{r}\right|_{r=r_{20}(t)}=-\hat{C}^{+}r_{20}'(t)-\frac{\hat{C}^{+}}{r_{20}(t)},
\end{align*}
where $\hat{C}^+=\int \dot{\omega}_{0}^{2}(z)dz$. In view of (\ref{SGT:1}) we obtain
$\varkappa_{1}=\varkappa_{2}=\hat{C}^{+}$.

Thus, in view of the corollary of Lemma~\ref{Lemma:2}, Section~5,
for the relation
$$
V^1_1\delta(x-\varphi_1)+V^2_1\delta(x-\varphi)
=\mathcal{O}_{\mathcal{D}'}(\varepsilon)
$$
to hold, it is sufficient to have
$V^1_1/r_{1}^{2}+V^2_2/r_{2}^{2}=0$ or, in more detail,
\begin{align}\label{VpV}
&(r_{1t}+ r_{2t})\left(B_\Omega+C_\Omega\right)
-\frac{2\beta_\tau \psi'_0}{\beta^2}B^z_\Omega\\
&+\left(\left.\frac{q}{r}\right|_{r=r_1}-\left.\frac{q}{r}\right|_{r=r_2}\right)C_\Omega
+\left(\frac{1}{r_{1}}+\frac{1}{r_{2}}\right)\left(\varkappa_{2}C_{\Omega}-\hat{C}\right)=0.\notag
\end{align}
Thus, we have proved the following assertion.

\begin{theorem}\label{theorem4}
The conditions of Theorem~{\rm1} and relations (\ref{ETA2}),
(\ref{BC}), (\ref{VpV}) are sufficient for the functions
$\check{\sigma}$, $\check{u}$ to be a weak asymptotic solution of
the phase field system~(\ref{PhHeat}), (\ref{PhAC}).
\end{theorem}

We consider relation (\ref{VpV}). It follows from the above that
$C_\Omega=C_\Omega(\eta)$ decreases sufficiently fast as
$\eta\to\infty$, in any case $\left| |\eta|
C_\Omega(\eta)\right|\leqslant\mathrm{const}$. This fact, the
estimate
$$
\hat{q}|_{r=r_1}-\hat{q}|_{r=r_2} =\mathcal{O}(|r_1-r_2|^\mu),\qquad \mu\in(0,1/2).
$$
proved in Lemma~\ref{lemma4} in Section~5, relation (\ref{52})
implies that the obvious estimate
$$
||r_1-r_2|^\mu C_\Omega| \leqslant
\varepsilon^{\mu}\mathrm{const}|\eta^{\mu} C_\Omega|
=\mathcal{O}(\varepsilon^\mu)
$$
for $t\geqslant t^*$.

For $t\leqslant t^{*}$, by Lemma~\ref{lemma6}, we have
\begin{equation*}
\left.q\right|_{r=r_{i}}-\left.q\right|_{r=r_{i0}}=
\mathcal{O}\left(|\psi_{0}\psi_{1}|^{\mu}\right),
\end{equation*}
where $\psi_{1}=r_{21}-r_{11}$, Hence, by (\ref{R1PH1}), and a statement similar to Lemma~\ref{Lemma:3}, in Section~5, we have
\begin{equation*}
\left|2(r_{10t}+r_{20t})-\left.\frac{q}{r}\right|_{r=r_{10}}+\left.\frac{q}{r}\right|_{r=r_{20}}+
\hat{C}^{+}\left(\frac{1}{r_{10}}+\frac{1}{r_{20}}\right)\right|=\mathcal{O}(\varepsilon^{\mu}),\quad t\leqslant t^{*}.
\end{equation*}
We introduce a function $V(\tau)\in C^{\infty}$, such that
$V'_{\tau}\in S(\mathbb{R}^{1})$, $V(-\infty)=0$, $V(\infty)=1$.
Then, in view of considerations similar to those used in Lemma~3, we
can show that the following estimate hold:
\begin{align*}
&\left(B_{\Omega}+C_{\Omega}\right)V(\tau)(r_{10t}+r_{20t})
+\left(-\left.\frac{q}{r}\right|_{r=r_{1}}+\left.\frac{q}{r}\right|_{r=r_{2}}\right)C_{\Omega}\\
&+\left(\varkappa_{2}C_{\Omega}-\hat{C}\right)
\left.\left(\frac{1}{r_{1}}+\frac{1}{r_{2}}\right)\right|=\mathcal{O}(\varepsilon^{\mu}).
\end{align*}

This implies that the left-hand side of (\ref{FCal:1}) is
estimated as $\mathcal{O}_{\mathcal{D}'}(\varepsilon^\mu)$ if
\begin{equation}\label{BpC}
\left(B_\Omega+C_\Omega\right)\left[(r_{10t}+r_{20t})
+r_{1t}+r_{2t}\right]
-2\frac{\beta'_\tau\psi'_{0}}{\beta^2}B^z_\Omega=0.
\end{equation}
From this relation we obtain
\begin{equation}\label{DifBpC}
\frac{\partial}{\partial\tau}(\tau(r_{11}+r_{21}))
=-(r_{10t}+r_{20t}) +2\frac{\beta'_\tau\psi'_{0}}{\beta^2}
\frac{B^z_\Omega}{B_\Omega+C_\Omega}.
\end{equation}
From this equation we determine the extensions of the functions
$r_{11}$ and $r_{21}$. We note that our argument results in an
equation that does not contain the temperature, namely, the
assumption that the classical solution exists till $t=t^*$ is
sufficient for constructing a global smooth approximation of the
solution.

Let us calculate the function $B_\Omega=B_\Omega(\eta)$
in more detail. We have
$$
B_\Omega(\eta)=\int\Omega'_z(\Omega'_z-\Omega'_\eta)\,dz.
$$
In this integral, we make the change of variable
$z\to-z-\eta$. Then we obtain
$$
\Omega'_\eta\to -(\Omega'_z-\Omega'_\eta), \qquad
(\Omega'_z-\Omega'_\eta)\to-\Omega'_z.
$$
Hence we have
$$
B_\Omega(\eta)=\int \Omega'_z \Omega'_\eta\,dz.
$$
Finally, we obtain
$$
B_\Omega=\frac12\int\big(\Omega'_z(\Omega'_z-\Omega'_\eta)
+\Omega'_z\Omega'_\eta\big)\,dz
=\frac12\int(\Omega'_z)^2\,dz.
$$
Thus, for finite $\eta$, the denominator $(B_\Omega+C_\Omega)$ in
the last term in the right-hand side of~(\ref{DifBpC}) does not
vanish. Moreover, using the explicit form of the functions
$\omega_0(z)$, $B^z_\Omega$, $B_\Omega$, and $C_\Omega$, we can
verify that
$$
2\frac{\beta'_\tau B^z_\Omega}{\beta^2}
=\mathcal{O}(|\eta|^{-N}),\qquad |\tau|\to\infty,
$$
where $N\gg1$ is an arbitrary number.

Therefore, we have
\begin{equation}\label{R1pR2}
r_{11}+r_{12}= -\left(r_{10t}+r_{20t}\right)\psi_{0}^{\prime-1} +\frac{2}{\tau}\int^{\tau}_{0} \frac{\beta'_\tau}{\beta^2}
\frac{B^z_\Omega}{B^z_\Omega+C_\Omega}\,d\tau'.
\end{equation}
We note that
\begin{equation}\label{DifR1pR2}
\psi'_{0}\left(1+\frac{\partial}{\partial\tau}
[\tau(r_{21}-r_{11})]\right)=\psi'_0\rho_\tau,
\end{equation}
where $\rho=\eta\beta^{-1}$ and the functions $\eta$ and $\beta$
are determined by Eqs.~(\ref{ETA2}) and~(\ref{BETA:1}).

The system of Eqs.~(\ref{R1pR2}) and (\ref{DifR1pR2}) allows one
to find the functions $r_{i1}$, $i=1,2$, and thus completely
determines the functions contained in the ansatz of the weak
asymptotic solution of system~(\ref{PhHeat}), (\ref{PhAC}).

It is easy to see that the solutions thus constructed satisfy
conditions~(\ref{BETAPHI:1}).

Next, using the explicit formulas for
$B^z_\Omega$, $B_\Omega$, and $C_\Omega$,
we can easily verify that $B^z_\Omega|_{\eta=0}=0$
(because $\omega_0(z)$ is odd) and
$(B_\Omega+C_\Omega)|_{\eta=0}\not=0$,
because $\dot\omega_0(z)$ is even.
This means that
\begin{equation}\label{63}
(r_{1t}+r_{2t})|_{\eta=0}=0.
\end{equation}

We note that the right-hand side of (\ref{Lq}) is a piecewise
smooth function for $|\psi|\geqslant\mathrm{const}>0$ and smooth
for $r\not=r_i$, $i=1,2$. As $\psi\to0$ ($t\to t^*$), the
right-hand side of (\ref{Lq}) becomes proportional to
$\delta(r-r^*_{0})$, where $r^*_0=r^*|_{\psi=0}$, and the
proportionality coefficient is negative. Thus, if $r$ and $t$ vary
in the respective neighborhoods of the points $r^*_0$ and $t^*$,
then the $\delta$-function with a negative coefficient appears and
disappears in the right-hand side of (\ref{Lq}), which results in
a negative soliton-like jump of the temperature in a neighborhood
of the point $t=t^*$, $r=r^*_0$.

{\it Calculation of the temperature jump}.
To prove the statement in the Introduction
concerning the temperature jump,
we must calculate the quantity
$$
[eI+\hat{q}+\check{q}]|_{t=t^*,\,r=r^*}.
$$

Let us verify that $\check{q}+\hat{q}$
is a continuous function.
This function is the sum of solutions
of the heat equation
with singular right-hand sides,
but these singularities arise
after the substitution of the continuous function
$e{\check{T}}$ (which, of course, is not a solution)
into the left-hand side.
It remains to prove that the solution differs
from the function $e{\check{T}}$ by a continuous function.
For this,
it suffices to verify that the right-hand sides,
which arise after the substitution of $e{\check{T}}$
into the equation,
do not generate any singularities
in the solution of the heat equation
in addition to those contained
in the function~$e{\check{T}}$.

Here we, in contrast to the preceding statements,
use the fact that the singularities
in the right-hand sides
arise as the result of the substitution.

By $\check{q}_i$, $i=1,2$, we denote the terms
in $\hat{q}+ \check{q}$ corresponding to the right-hand sides:
\begin{align*}
f_1&=-\left(\frac{\partial^2}{\partial x^2}\gamma^-
\frac{(\varphi_1-x)(x-\varphi_2)}{\psi}\right)
B(\tau)[H(r-r_1)-H(r-r_2)],
\\
f_2&=-\left(\frac{\partial^2}{\partial x^2}\hat{\gamma}
\frac{(\varphi_1-x)(x-\varphi_2)}{\psi}\right)
(1-B(\tau))[H(r-r_2)-H(r-r_1)],
\end{align*}
The other terms in $q$ are solutions of the heat
equation with piecewise continuous right-hand side and hence
are continuous.

The functions $\hat{q}_i$ have a similar property if
they are calculated up to $\mathcal{O}(\varepsilon)$.
Indeed, we denote
$$
\Pi=\gamma^-\frac{(r_1-r)(r-r_2)}{\psi}
$$
and represent, for example,
the function $\hat{q}_1$ in the form
$$
\hat{q}_1=-\frac1{2\sqrt{2\pi}}
\int^t_0\frac{B(\tau(t',\varepsilon))}{\sqrt{t-t'}}
\int^{\varphi_2}_{\varphi_1}\frac{\partial^2\Pi}{\partial \xi^2}
e^{-\frac{(x-\xi)^2}{4(t-t')} }\,d\xi dt'.
$$

Omitting the number factor and integrating by parts
in the integral over~$\xi$, we obtain
\begin{align*}
\hat{q}_1&=-\int^t_0\left.\frac{B}{\sqrt{t-t'}}
e^{-\frac{(r-\xi)^2}{4(t-t')}} \right|^{r_2}_{r_1}
-\int^t_0\frac{B}{\sqrt{t-t'}}\frac{\partial^2}{\partial \xi^2}
e^{-\frac{(r-\xi)^2}{4(t-t')}}\,d\xi dt'
\\
&=-\int^t_0\frac{B}{\sqrt{t-t'}}
\left[e^{-\frac{(r-\xi)^2}{4(t-t')}}\gamma^-_2
-e^{-\frac{(r-\xi)^2}{4(t-t')}}\gamma^+_1\right]\,dt'
\\
&\qquad +\int^t_0\frac{\partial}{\partial t'}
\left(\int^{r_2}_{r_1} \frac{\Pi
e^{-\frac{(r-\xi)^2}{4(t-t')}}}{\sqrt{t-t'}}\,d\xi\right)\,dt'
+\int^t_0\left(\int^{r_2}_{r_1} B \frac{\partial \Pi}{\partial t'}
\frac{e^{-\frac{(r-\xi)^2}{4(t-t')}}}{\sqrt{t-t'}}\,d\xi\right)\,dt'.
\end{align*}

The last term is a solution of the heat equation
with the piecewise continuous right-hand side
$$
B\frac{\partial \Pi}{\partial t} [H(r-r_1)-H(r-r_2)],
$$
and hence it is continuous
uniformly w.r.t $\varepsilon\geq0$.
The other terms are also continous functions
uniformly w.r.t $\varepsilon\geq0$.

It is easy to see that in the formula for $\hat{q}_1$
there is no term containing the derivative
$\frac{\partial B}{\partial t}$,
because this term would be of order $\mathcal{O}(\varepsilon)$.

Indeed
$$
\left|\frac{\partial B}{\partial t}\Pi [H(r-r_1)-H(r-r_2)]\right|
\leqslant C|B'\cdot\psi'_{0}\cdot\rho|,
$$
since $|\Pi|\leqslant C\psi$ for $r\in[r_1,r_2]$ and
$\varepsilon^{-1}\psi=\rho$. The derivative $B'_\rho$ decreases
faster than any power of $|\rho|^{-1}$,
$|\rho'_\tau|\leq\mathrm{const}$. This implies that
$$
\int^t_0\frac{\partial B}{\partial t}\,dt'
=\mathcal{O}(\varepsilon).
$$

Hence we have
$$
[\check{q}+\hat{q}]|_{r=r^*,\,t=t^*}=0.
$$
Let us calculate the function $eI$.
We have
$$
eI|_{r=r^*}=
\frac{1}{4}\left(\left(r_{1}^2\right)_{t}-\left(r_{2}^2\right)_{t}\right)
=-\frac{\psi_{0}'\dot{\rho}_{\tau}}{4}\left(r_{1t}+r_{2t}\right).
$$
It follows from (\ref{RhoSt:1}) that
\begin{gather*}
\dot\rho_\tau\to1,\qquad \tau\to\infty,
\\
\dot\rho_\tau\to0,\qquad \tau\to-\infty
\end{gather*}
and hence
$$
[eI]|_{\substack{r=r^*,\\t=t^*}} =-\frac14\lim_{t\to
t^*-0}(r_{10t}^{2}+r_{20t}^{2}).
$$

So, the jump of the temperature has the form
$$
[\bar{\theta}]|_{\substack{r=r^*,\\t=t^*}} =-\frac{r_{10}(t^*)}4\lim_{t\to t^*-0}(r_{10t}^{2}+r_{20t}^{2}).
$$

\section{Technique of the weak asymptotics\\ method}\label{WAsMeth}

First, we recall the definition of
\emph{regularization of the generalized function}.

\begin{definition}\label{Def:Regul}
A family of functions $f(x,\varepsilon)$
smooth for $\varepsilon>0$ and satisfying the condition
\begin{equation*}
\underset{\varepsilon\to0}{\mathrm{w}-\lim}f(x,\varepsilon)=f(x)
\end{equation*}
is called the regularization of the generalized function$f(x)$.
\end{definition}

We note that, by definition,
the last relation can be rewritten as
\begin{equation*}
\lim_{\varepsilon\to0}\langle
f(x,\varepsilon),\zeta(x)\rangle=\langle f,\zeta\rangle
\end{equation*}
for any test function $\zeta(x)$
(from now on, $\langle\ ,\ \rangle$
denotes the action of a generalized function
on a test function).

\begin{lemma}
Let $\Gamma_{t}=\{x-\varphi(t)=0\}$,
$x\in\mathbb{R}$, where $\varphi(t)$ is a smooth function,
let $\omega(z)\in\mathbb{S}$
($\mathbb{S}$ is the Schwartz space),
and let $\beta=\beta(t)>0$.
Then the following relation holds for any test function
$\zeta(x)${\rm:}
$$
\frac{1}{\varepsilon}
\Big\langle\omega\left(\beta\frac{x-\varphi(t)}{\varepsilon}\right),
\zeta(x)\Big\rangle
=\frac{1}{\beta}A_{\omega}\zeta(\varphi)+\mathcal{O}(\varepsilon),
$$
where $A_{\omega}=\int^{\infty}_{-\infty}\omega(z)dz$.
\end{lemma}

\begin{proof}
The expression in the right-hand side can be written as
\begin{equation*}
\frac{1}{\varepsilon}
\int_{\mathbb{R}^{1}}\omega
\left(\beta\frac{x-\varphi(t)}{\varepsilon}\right)\zeta(x)\,dx
= \frac{\zeta(\varphi)}{\beta}\int_{\mathbb{R}^{1}}\omega(z)\,dz
+\mathcal{O}(\varepsilon).
\end{equation*}
Here we perform the change of variables
$z=\beta(x-\varphi)/\varepsilon$ and apply the Taylor formula
to the integrand at the point $x=\varphi$.
By definition, the last integral is the action of the
generalized function
$\beta^{-1}A_{\omega}\delta(x-\varphi)$ on the test function $\zeta$.
\end{proof}

\textbf{Proof of the formula in Example~1 in the
Introduction.}
Let $\omega(z) \in C^{\infty}$,
$\omega'\in\mathbb{S}(\mathbb{R}^{1})$,
$\lim_{z\to+\infty}\omega(z)=0$,
and
$\lim_{z\to-\infty}\omega(z)=1$.
We verify that
\begin{equation}\label{HW:1}
\omega\left(\frac{x-x_{0}}{\varepsilon}\right)-H(x-x_{0})
=\mathcal{O}_{\mathcal{D}'}(\varepsilon),
\quad x_{0}=\mathrm{const},
\end{equation}
where $H$ is the Heaviside function.
By Definition~\ref{WkMDef},
we consider the expression
\begin{align}\label{HW:2}
&\int_{\mathbb{R}^{1}}
\left[\omega\left(\frac{x-x_{0}}{\varepsilon}\right)
-H(x-x_{0})\right]\zeta(x)\,dx\\
&\qquad
=\varepsilon\zeta(x_{0})
\int_{\mathbb{R}^{1}}\left[\omega(z)-H(z)\right]\,dz
+\mathcal{O}(\varepsilon^{2}).\notag
\end{align}
Here we performed the change of variables
$z=(x-x_{0})/\varepsilon$ and applied the Taylor formula
to the functions $\zeta(x)$
at the point $x=x_{0}$.
In view of our assumptions and the properties of the
function~$H$, the integral in the right-hand side of the last
relation converges and hence the right-hand side of (\ref{HW:2})
is of order  $\mathcal{O}(\varepsilon)$.
We thus obtain estimate (\ref{HW:1}).
\medskip

\textbf{Proof of the formula in Example~2 in the
Introduction.}
Let $\omega_{i}(z) \in C^{\infty}$,
$\left(\omega_{i}\right)'\in\mathbb{S}(\mathbb{R}^{1})$,
$\lim_{z\to+\infty}\omega_{i}(z)=0$,
and $\lim_{z\to-\infty}\omega_{i}(z)=1$, $i=1,2$.
We consider the integral
\begin{align*}
J&\buildrel{\rm def}\over=
\int_{\mathbb{R}^{1}}
\omega_{1}\left(\frac{x-x_{1}}{\varepsilon}\right)
\omega_{2}\left(\frac{x-x_{2}}{\varepsilon}\right)\zeta(x)\,dx\\
&=\int_{\mathbb{R}^{1}}
\omega_{1}\left(\frac{x-x_{1}}{\varepsilon}\right)
\omega_{2}\left(\frac{x-x_{2}}{\varepsilon}\right)
\left(\int_{-\infty}^{x}\zeta(y)dy\right)'_x\,dx.
\end{align*}
Integrating by parts, we obtain
\begin{align*}
J&=-\int_{\mathbb{R}^{1}}
\dot{\omega}_{1}\left(\frac{x-x_{1}}{\varepsilon}\right)
\omega_{2}\left(\frac{x-x_{2}}{\varepsilon}\right)
\left(\int_{-\infty}^{x}\zeta(y)dy\right)\,dx\\
&\qquad
-\int_{\mathbb{R}^{1}}
\dot{\omega}_{2}\left(\frac{x-x_{2}}{\varepsilon}\right)
\omega_{1}\left(\frac{x-x_{1}}{\varepsilon}\right)
\left(\int_{-\infty}^{x}\zeta(y)dy\right)\,dx.
\end{align*}
We perform the change of variables $z=(x-x_{i})/\varepsilon$
with $i=1$ in the first integral
and with $i=2$ in the second integral
and apply the Taylor formula to the function
$F(x)=\int^{x}_{-\infty}\zeta(y)\,dy$
at the points $x=x_{i}$, $i=1,2$, respectively.
Then we calculate the first and the second integral,
we have
\begin{gather*}
J=B_{1}\left(\frac{\triangle
x}{\varepsilon}\right)H(x-x_{1})+B_{2}\left(\frac{\triangle x}{\varepsilon}\right)
H(x-x_{2})+\mathcal{O}_{\mathcal{D}'}(\varepsilon),\quad
\triangle x=x_{1}-x_{2},
\\
B_{1}(\rho)=\int_{\mathbb{R}^{1}}
\dot{\omega}_{1}(z)\omega_{1}(-z-\rho)\,dz,\quad
B_{2}(\rho)=\int_{\mathbb{R}^{1}}
\dot{\omega}_{2}(z)\omega_{1}(z-\rho)\,dz.
\end{gather*}

To calculate the linear combination of generalized functions
up to $\mathcal{O}_{\mathcal{D}'}(\varepsilon^{\alpha})$,
we must improve the classical definition of linear independence.
This improvement plays the key role in the study
of problems with interaction of nonlinear waves.

Indeed, let $\phi_{1}\neq\phi_{2}$ be independent of~$x$.
We consider the expression
\begin{equation}\label{LinDelta:1}
g_{1}\delta(x-\phi_{1})+g_{2}\delta(x-\phi_{2})
=\mathcal{O}_{\mathcal{D}'}(\varepsilon^{\alpha}),\qquad
\alpha>0,
\end{equation}
where the functions $g_{i}$ are independent of $\varepsilon$.
Clearly, for the last relation to be satisfied,
it suffices to have
\begin{equation*}
g_{i}=\mathcal{O}(\varepsilon^{\alpha}),\ \ \ i=1,2,
\end{equation*}
or, with the properties of the functions $g_{i}$ taken into account,
\begin{equation*}
g_{i}=0,\qquad i=1,2.
\end{equation*}

But, if we assume that the coefficients $g_{i}$
depend on the parameter $\varepsilon$,
then the above estimates do not work.
Namely, let us consider the following specific case
of this dependence:
\begin{equation}\label{g:g}
g_{i}=\phi_{i}^{n}S_{i}(\triangle\phi/\varepsilon),\quad n>1,\quad
i=1,2,
\end{equation}
where the functions $S_{i}(\sigma)$ decrease sufficiently fast as
$|\sigma|\to\infty$, and $\triangle\phi=\phi_{2}-\phi_{1}$.

\begin{lemma}[\rm Linear independence
of generalized functions]\label{Lemma:2}
Suppose that the estimate
\begin{equation*}
|\sigma S_{i}(\sigma)|\leqslant \mathrm{const},\quad i=1,2, \quad
-\infty<\sigma<+\infty
\end{equation*}
holds. Then, for $\alpha=1$, expression {\rm(\ref{LinDelta:1})}
implies the relation
\begin{equation}\label{LinCoef:1}
S_{1}=-S_{2}.
\end{equation}
\end{lemma}

\begin{proof}
Using the Taylor formula in (\ref{LinDelta:1})
and taking (\ref{g:g}) into account, we obtain
\begin{align*}
\phi_{1}^{n}S_{1}\zeta(\phi_{1})+\phi_{2}^{n}S_{2}\zeta(\phi_{2})
=&\phi_{1}^{n}S_{1}\zeta(\phi_{1})+\phi_{2}^{n}S_{2}\zeta(\phi_{1})\\
&+\phi_{2}^{n}S_{2}(\phi_{2}-\phi_{1})\zeta'(\phi_{1}+\mu\phi_{2}),
\end{align*}
where $0<\mu<1$
Since the function $\sigma S_{2}(\sigma)$ is uniformly bounded
in $\sigma\in\mathbb{R}^{1}$, we obtain
\begin{equation*}
S_{2}(\triangle\phi/\varepsilon)(\phi_{2}-\phi_{1})
=\left.\{-\sigma S_{2}(\sigma)\}
\right|_{\sigma=\triangle\phi/\varepsilon}\cdot\varepsilon
=\mathcal{O}(\varepsilon).
\end{equation*}
If $S_{1}=-S_{2}$ then expression (\ref{LinDelta:1}) can be
rewritten as
\begin{equation*}
(\phi_{1}^{n}S_{1}+\phi_{2}^{n}S_{2})\zeta(\phi_{1})
=(\phi_{1}^{n}-\phi_{2}^{n})S_{1}=\mathcal{O}(\varepsilon).
\end{equation*}
Thus we obtain the assertion of  Lemma~\ref{Lemma:2}.
\end{proof}

\begin{corollary}
Suppose that
$$
|\sigma S_i(\sigma)|\leq\mathrm{const},\qquad
i=1,2,\quad\sigma\geq0
$$
and the functions $\phi_i=\phi_i(t,\varepsilon)$
are continuous and uniformly continuous and satisfy
the condition that
$\Delta\phi>0$ for $t>0$ and $\varepsilon\geq0$
and $\Delta\phi=\mathcal{O}(\varepsilon)$ for $t<0$.
Then, for $\alpha=1$, expression {\rm(66)}
implies~{\rm(64)}.
\end{corollary}

\begin{proof}
Following the above argument, we obtain
$$
S_1\zeta(\phi_1)+S_2\zeta(\phi_2)
=(\phi_{1}^{n}-\phi_{2}^{n})S_{1}
+\phi_{2}^{n}S_2(\phi_2-\phi_1)\zeta'(\phi_1+\mu\phi_2).
$$

In view of our assumptions,
$$
\phi_{2}^{n}S_2(\phi_2-\phi_1)=\mathcal{O}(\varepsilon),\quad
(\phi_{1}^{n}-\phi_{2}^{n})S_{1}=\mathcal{O}(\varepsilon)
$$
uniformly in~$t$, which implies the desired assertion.
\end{proof}

\begin{lemma}\label{Lemma:3}
Suppose that $f(t)\in\mathbb{C}^{1}$, $f(t_{0})=0$,
and $f'(t_{0})\neq 0$.
Suppose also that $g(t,\tau)$ locally uniformly in $t$
satisfies the condition
\begin{equation*}
|\tau g(t,\tau)|\leqslant \mathrm{const},\quad
|\tau g'_t(t,\tau)|\leqslant \mathrm{const},\quad -\infty<\tau<\infty,
\end{equation*}
and $g(t_{0},\tau)=0$. Then the inequality
\begin{equation*}
\left|g\left(t,f(t)/\varepsilon\right)\right|\leqslant \varepsilon
C_{\hat{t}},
\end{equation*}
where $C_{\hat{t}}=\mathrm{const}$,
holds in any interval $0\leqslant t\leqslant \hat{t}$ that does
not contain zeros of the function $f(t)$ except for~$t_{0}$.
\end{lemma}

\begin{proof}
The fraction $\frac{f(t)}{t-t_{0}}$ is locally bounded in~$t$.
The fraction $\frac{\tau g(t,\tau)}{t-t_{0}}$
is also locally bounded.
We have
\begin{equation*}
g\left(t,f(t)/\varepsilon\right)=
\varepsilon\cdot\frac{g\left(t,f(t)/\varepsilon\right)}{(t-t_{0})}\cdot
\frac{f(t)}{\varepsilon}\cdot\frac{t-t_{0}}{f(t)}.
\end{equation*}
According to the assumptions of the lemma,
the last factor in the right-hand side is bounded
on the interval under study.
The product of the first and second factors
(without $\varepsilon$) is bounded in view of the properties
of the function $g(t,\tau)$.
\end{proof}

\begin{corollary}\label{Corollary:2}
Suppose that the assumptions of Lemma~{\rm\ref{Lemma:3}}
are satisfied for $0\leqslant\tau<\infty$
$(-\infty<\tau\leqslant0)$.
Then the assertion of Lemma~{\rm\ref{Lemma:3}}
holds on any half-interval $(t_{0},\hat{t}]$,
which does not contain zeros of the function~$f(t)$,
and $\mathrm{sign} \hat{t}=\mathrm{sign} f(t)$,
$t\in (t_{0},\hat{t}]$.
\end{corollary}

The proof of Corollary~\ref{Corollary:2} is obvious.

\begin{lemma}\label{lemma4}
The following inequality holds{\rm:}
$$
\hat{q}|_{r=r_1}-\hat{q}|_{r=r_2} =\mathcal{O}(|r_1-r_2|^\mu).
$$
\end{lemma}

\begin{proof}
We consider only one of the the functions $\hat{q}$ from
(\ref{hatq}), namely, the functions determined by relation
(\ref{31}), and apply the method developed in \cite{14}. We
consider the difference of expressions (\ref{hatq}) omitting the
number factors:
\begin{align*}
&\hat{q}|_{r=r_1}-\hat{q}|_{r=r_2}
\\
&\qquad=\int^t_0 \int^{r_2}_{ri_1} \frac{\displaystyle
[\exp\frac{(r_1-\xi)^2}{4(t-t')}
-\exp(-\frac{(r_2-\xi)^2}{4(t-t')})]}
{\psi\sqrt{t-\tau}}\,d\xi\,d\tau
\\
&\qquad=\int^t_0 \int^{r_2}_{\frac{r_1+r_2}2} \frac{\displaystyle
\exp(-\frac{(r_2-\xi)^2}{4(t-t')})
[\exp\frac{(r_2-r_1)(r_1+r_2-2\xi)}{4(t-t')} -1]}
{\psi\sqrt{t-t'}}
\\
&\qquad\qquad\times \frac{(t-t')^\mu}{(r_2-r_1)^\mu}
\frac{(r_2-r_1)^\mu}{(t-t')^\mu} \,d\xi\,d\tau
\\
&\qquad=\int^t_0 \int_{r_1}^{\frac{r_1+r_2}2} \frac{\displaystyle
\exp(-\frac{(r_1-\xi)^2}{4(t-t')})
[\exp\frac{(r_1-r_2)(r_1+r_2-\xi)}{4(t-t')} -1]} {\psi\sqrt{t-t'}}
\\
&\qquad\qquad\times \frac{(t-t')^\mu}{(r_2-r_1)^\mu}
\frac{(r_2-r_1)^\mu}{(t-t')^\mu} \,d\xi\,dt',
\end{align*}
where we choose $\mu\in(0,1/2)$.

Next, using the inequality
$$
\bigg|\frac{e^{-\alpha r}-1}{r^\mu}\bigg|\leqslant\mathrm{const}
$$
for $\alpha>0$, $r\in[0,\infty)$, and taking into account the fact
that the integral
$$
\int^t_0 (t-t')^{-(1/2+\mu)}\,dt',\qquad
\mu\in(0,1/2),
$$
converges, we obtain the statement of the lemma for $r_2\geqslant
r_1$. If $r_2\leqslant r_1$ (this is true for $t\geq t^*$), then
we must change the exponents outside the square brackets in the
integrands.
\end{proof}

\begin{lemma}\label{lemma5}
Suppose that $\omega(z)\in C^\infty$
decreases faster than any power of $|z|^{-1}$ as $|z|\to\infty$.
Then
$$
\varepsilon^{-1}\omega((r-r_i)/\varepsilon) \hat{q}_j
=\hat{q}_j\delta(r-r_i)\int\omega(z)\,dz
+\mathcal{O}_{\mathcal{D}'}(\varepsilon^\mu),
$$
where $\hat{q}_j$ is one of the functions defined by the relations
(\ref{qqqq}), $i=1,2$, $\mu\in(0,1/2)$.
\end{lemma}

\begin{proof}
We consider one of the functions $\hat{q}$, namely,
the function defined by the relations
$$
\hat{q}=\frac1{2\sqrt{\pi}}\int^t_0\frac{dt'}{\psi\sqrt{t-t'}}
\int^{r_2}_{r_1}\exp(-\frac{(r-\xi)^2}{4(t-t')})\,d\xi.
$$
The desired relation can be written as
$$
\varepsilon^{-1}\int \zeta(r)
\omega\bigg(\frac{r-r_i}{\varepsilon}\bigg)\hat{q} dr
=\zeta(r_i)\hat{q}|_{r=r_i} \int\omega(z) dz
+\mathcal{O}(\varepsilon^\mu)
$$
or, omitting the number factors, in the form
\begin{align}\label{71}
I&\od \varepsilon^{-1}\int\zeta(r)
\omega\bigg(\frac{(r-r_i)}{\varepsilon}\bigg)\hat{q} dr
\\
&\qquad\times \int^t_0 \frac{dt'}{\psi\sqrt{t-t'}}
\int^{r_2}_{r_1} \exp\bigg(-\frac{(r-\xi)^2}{4(t-t')}\bigg)\,d\xi\
dr
\notag\\
&=\zeta(r_i)\int^t_0 \frac{dt'}{\psi\sqrt{t-t'}} \int^{r_2}_{r_1}
\exp\bigg(-\frac{(r_i-\xi)^2}{4(t-t')}\bigg) d\xi
+\mathcal{O}(\varepsilon^\mu). \notag
\end{align}
In the left-hand side, we change the variables
$(r-r_i)/\varepsilon=z$ and obtain
$$
I=\int\zeta(r_i+\varepsilon z)\omega(z) \int^t_0
\frac{dt'}{\psi\sqrt{t-t'}} \int^{r_2}_{r_1}
\exp\bigg(-\frac{(r_i-\xi+\varepsilon z)^2}{4(t-t')}\bigg) dz
d\xi.
$$
It is clear that to prove relation (67),
it suffices to prove that, with accuracy up to small values,
the term $\varepsilon z$ can be omitted in the exponent.

We consider the expression
\begin{align*}
J&=\int\zeta(r_i+\varepsilon z)\omega(z)
\int^t_0\frac{dt'}{\psi\sqrt{t-t'}}
\\
&\qquad\times \int^{r_2}_{r_1}
\bigg[\exp\bigg(-\frac{(r_i-\xi+\varepsilon z)^2}{4(t-t')}\bigg)
-\exp\bigg(-\frac{(r_i-\xi)^2}{4(t-t')}\bigg)\bigg] dz d\xi.
\end{align*}
Since the function $\omega$ decreases fast,  we can assume
that $|z|<\varepsilon^{-\delta}$, $\delta\in(0,1)$.
Next it suffices to consider the integral over $t'$
from $0$ to $t-\varepsilon^{1-\delta}$, because if
$t'\in[t-\varepsilon^{1-\delta},t]$, then the obtained integral
can be estimated as $\mathcal{O}(\varepsilon^{(1-\delta)/2})$.

In the remaining integrals, we can use the same method as in the
proof of Lemma~\ref{lemma4}. Namely, we transform the difference
of the exponential functions:
\begin{align*}
&\exp\bigg(-\frac{(r_i-\xi+\varepsilon z)^2}{4(t-t')}\bigg)
-\exp\bigg(-\frac{(r_i-\xi)^2}{4(t-t')}\bigg)
\\
&\qquad =\bigg[\exp\bigg(-\frac{(2(r_i-\xi)-\varepsilon
z)\varepsilon z} {4(t-t')}\bigg)-1 \bigg] \bigg(\frac{\varepsilon
z}{t-t'}\bigg)^\gamma \bigg(\frac{\varepsilon
z}{t-t'}\bigg)^{-\gamma},
\\
&\qquad
\gamma\in(0,1/2),
\end{align*}
and apply the estimate
$$
\bigg|\bigg(\exp\bigg( -\frac{(2(r_i-\xi)-\varepsilon
z)\varepsilon z} {4(t-t')}\bigg)-1\bigg)\bigg/
\bigg(\frac{t-t'}{\varepsilon z}\bigg)^\gamma\bigg|<1,
$$
which, obviously, is valid in the required range of variables.
Hence, as in Lemma~4, we obtain the estimate
$$
J=\mathcal{O}(\varepsilon^\mu),\qquad \mu\in(0,1/2),
$$
which proves relation~(\ref{71}).
\end{proof}

Similarly to Lemmas~\ref{lemma4} and~\ref{lemma5}, we prove the
following statement.

\begin{lemma}\label{lemma6}
The following relation holds for any of the functions $\hat q_i$
defined by relations (\ref{qqqq}):
$$
\hat q_i(r_i,t)-\hat q_i(\hat q_{i0},t)
=\mathcal{O}(|\psi_0\psi_1|^\mu),\qquad \mu\in(0,1/2).
$$
\end{lemma}

\begin{lemma}\label{lemma7}
Relations (\ref{VVV}), (\ref{BVVV}) hold.
\end{lemma}

\begin{proof}
In Section~3, we introduced the function
$\Omega(z,\eta)$.
Of course, we can change the signs in the arguments
using the fact that the function $\omega_0(z)$ is odd,
but the form presented above is convenient for calculations.

Calculating the derivatives $u_t$ and $u_x$, we obtain
\begin{align}\label{UdifT:1}
&\check{u}_{t}=\frac{1}{2\varepsilon}
\left\{\left[\beta(r_{1}-r)\right]_{t}
\dot{\omega}_{0}\left(\beta\frac{r_{1}-r}{\varepsilon}\right)
+\left[\beta(r-r_{2})\right]_{t}
\dot{\omega}_{0}\left(\beta\frac{r-r_{2}}{\varepsilon}\right)
\right.\notag\\
&\qquad-\left[\beta(r_{1}-r)\right]_{t}
\dot{\omega}_{0}\left(\beta\frac{r_{1}-r}{\varepsilon}\right)
\omega_{0}\left(\beta\frac{r-r_{2}}{\varepsilon}\right)\\
&\qquad-\left.\left[\beta(r-r_{2})\right]_{t}
\dot{\omega}_{0}\left(\beta\frac{r-r_{2}}{\varepsilon}\right)
\omega_{0}\left(\beta\frac{r_{1}-r}{\varepsilon}\right)\right\}\notag,
\end{align}
\begin{align}\label{UdifX:1}
&\check{u}_{r}=\frac{\beta}{2\varepsilon}\left\{ -\dot{\omega}_{0}
\left(\beta\frac{r_{1}-r}{\varepsilon}\right)+
\dot{\omega}_{0}\left(\beta\frac{r-r_{2}}{\varepsilon}\right)
\right.\\
&\qquad\left.+\dot{\omega}_{0}
\left(\beta\frac{r_{1}-r}{\varepsilon}\right)
\omega_{0}\left(\beta\frac{r-r_{2}}{\varepsilon}\right)-
\dot{\omega}_{0}\left(\beta\frac{r-r_{2}}{\varepsilon}\right)
\omega_{0}\left(\beta\frac{r_{1}-r}{\varepsilon}\right)\right\}
\notag.
\end{align}
Now we group the terms in the subintegral expression so that one
group contain terms with the factor
$\dot\omega_0(\beta(r_i-r)/\varepsilon)$, and the other group
contain terms with the factor
$\dot\omega_0(\beta(r-r_2)/\varepsilon)$. Thus, we write the
integral as the sum of two integrals. Then we perform the change
of variables
\begin{equation}
r\to r_1-\frac{\varepsilon z}{\beta}
\end{equation}
in the first integral
and the change of variables
\begin{equation}\label{rtz}
r\to r_2+\frac{\varepsilon z}{\beta}
\end{equation}
in the second integral.

From (\ref{UdifT:1}) and (\ref{UdifX:1}), we obtain
\begin{align}\label{uuuu}
&\int r^{2} u_t u_r\zeta(r) dr
\\
&\quad
=-\frac{r^{2}}{2}\bigg[r_{1t}+\frac{\beta_\tau}{\beta^2}\psi'_0\bigg]
\zeta(r_1)[(1+z)(\Omega'_{z}-\Omega'_\eta)
+(1-z-\eta)\Omega'_{\eta}][\Omega'_\eta-\Omega'_z]\,dz
\notag\\
&\quad
-\frac{r^{2}}{2}\bigg[r_{2t}+\frac{\beta_\tau}{\beta^2}\psi'_0\bigg]
\zeta(r_2)[(1+z)\Omega'_{1z}
+r^{2}(1-z-\eta)\Omega'_{\eta}][\Omega'_\eta-\Omega'_z]\,dz
\notag\\
&\quad
+\mathcal{O}(\varepsilon).
\notag
\end{align}

Here we used the formula for the derivative
$\beta_t=\varepsilon^{-1}\beta_\tau\psi'_{0t}$ and the fact that
the functions in the integrand in~(\ref{uuuu}) are invariant under
the change $z\to-z-\eta$, since each bracket is multiplied by~$-1$
in this change.

Now we calculate the other expressions contained in
(\ref{Defin:2}) have
\begin{align}\label{urz}
&\frac{\varepsilon}{2}\int r^{2} \zeta'(r)(\check u_r)^2 dr
=\frac{\beta^2}{8\varepsilon}\int r^{2}\zeta'(r)
\bigg[-\dot{\omega}_{0}
\left(\beta\frac{r_{1}-r}{\varepsilon}\right)
+\dot{\omega}_{0}\left(\beta\frac{r-r_{2}}{\varepsilon}\right)
\\
&\qquad + \dot{\omega}_{0}
\left(\beta\frac{r_{1}-r}{\varepsilon}\right)
\omega_{0}\left(\beta\frac{r-r_{2}}{\varepsilon}\right) -
\dot{\omega}_{0}\left(\beta\frac{r-r_{2}}{\varepsilon}\right)
\omega_{0}\left(\beta\frac{r_{1}-r}{\varepsilon}\right)\bigg]^2 dr
\notag\\
&=\frac{\beta}{8}\big(r_{1}^{2}\zeta'(r_1)+r_{2}^{2}\zeta'(r_2)\big)
\int\big\{\dot\omega_0(z)(1-\omega_0(-z-\eta))^2
\notag\\
&\qquad
-\dot\omega_0(z)\dot\omega_0(-z-\eta)
(1-\omega_0(z))(1-\omega_0(-z-\eta)) \big\}\,dz
+\mathcal{O}(\varepsilon).
\notag
\end{align}
It is easy to see that the integrand in the right-hand side of
(\ref{urz}) can be written as
\begin{align*}
&\int\{(\dot\omega_0(z)(1-\omega_0(-z-\eta))^2
\\
&\qquad
-\dot\omega_0(z)\dot\omega(-z-\eta)
(1-\omega_0(z))(1-\omega_0(-z-\eta))\}\,dz
=\frac14\int(\Omega'_z)^2\,dz.
\end{align*}
Thus, we finally obtain
$$
\frac{\varepsilon}{2}\int r^{2}\zeta'(\check{u}_r)^2\,dx
=\frac{\beta}{4}(r_{1}^{2}\zeta'(r_1)+r_{2}^{2}\zeta'(r_2))
\int(\Omega'_z)^2\,dz+\mathcal{O}(\varepsilon).
$$
Then the left-hand side of (\ref{urz}) is positive, and hence the
integrand expression in the right-hand side is also positive
(until the measure of the points at which
$|u_r|\geqslant\mbox{const}$ is positive).

Similarly we obtain:
$$
\frac1{\varepsilon}\int r^{2}F(\check u)\zeta'_r dr
=\frac{\beta^{-1}}{2}(r_{1}^{2}\zeta'(r_1)+r_{2}^{2}\zeta'(r_2))
\int F(\Omega) dz +\mathcal{O}(\varepsilon).
$$

To prove this relation, we must successively perform the change of
variables (\ref{rtz}) in the integral in the left-hand side and
then consider the half-sum of the integrals obtained, which,
obviously, is equal to the original integral in the left-hand
side. Moreover, we must take into account the estimates
(\ref{r34}), (\ref{r35}), and
$$
\Omega(z,\eta)
=\begin{cases}
1+\mathcal{O}(\exp(-2z)), &z\to\infty,
\\
\omega_0(-z-\eta)+\mathcal{O}(\exp(2z)), &z\to-\infty,
\end{cases}
$$
which together with the explicit form of the function
$F(u)=\frac{u^4}{4}-\frac{u^2}{2}+\frac14$
imply that the integrals
$$
\int zF(\Omega)\,dz
$$
converge.
\end{proof}


\begin{thebibliography}{99}

\bibitem{1}
\textsc{X.Chen} \emph{Spectrum for the Allen-Cahn, Cahn-Hilliard
and Phase-field equations for generic interfaces.} 1994, Commun.
in Part. Dif. Equat. 19(7), pp. 1371--1395.

\bibitem{12}
\textsc{G.Caginalp} \emph{Stefan and Hele-Shaw type models as
asymptotic limits of the phase-field equations.} Phys. Rev. Vol.
39, pp. 5887-5896, 1989.

\bibitem{9}
\textsc{V.G.Danilov} \emph{Propagation and interaction of shock
waves of quasilinear equation.} Nonlinear Studies~8, 2001, no~1,
211-245.

\bibitem{6}
\textsc{V.G.Danilov} \emph{Generalized solutions describing
singularity interaction.} 2002, Int. J. Math. and Math. Sci. 29,
no. 8, 481--494.

\bibitem{2}
\textsc{V.G.Danilov, G.A.Omel'yanov \& E.V.Radkevich}
\emph{Asymptotic solution of a phase field system and the modified
Stefan problem.} 1995, Differential'nie Uravneniya 31(3), 483--491
(in Russian).(English translation in Differential Equations 31(3),
1993.

\bibitem{4}
\textsc{V.G.Danilov, G.A.Omel'yanov \& E.V.Radkevich}
\emph{Hugoniot-type conditions and weak solutions to the
phase-field system.} 1999, Euro Journal of Applied Mathematics,
vol.~10, pp.~55--77.

\bibitem{10}
\textsc{V.G.Danilov and V.M.Shelkovich} \emph{Propagation and
interaction of nonlinear waves.} Eighth International Conference
on Hyperbolic Problems, Theory-Numerics-Applications, Abstracts,
Magdeburg, Germany, February 28-March 3, 2000, pp. 326-328.

\bibitem{11}
\textsc{V.G.Danilov G.A.Omel'yanov and V.M.Shelkovich}
\emph{Weak asy\-mptotics method and interaction of nonlinear
waves.} Amer. Math. Soc. Transl. (2) Vol. 208, 2003.

\bibitem{5}
\textsc{A.M.Meirmanov} \emph{The Stefan problem.} 1986, "Nauka",
Novosibirsk; English transl., de Gruyter, Berlin, 1992.

\bibitem{111}
\textsc{A.Meirmanov} \emph{The Stefan problem with surface
tension in the three-dimensional case with spherical symmetry:
nonexistence of the classical solution.}
Euro. J. Appl. Math.,  Vol.~5, pp. 1--19, 1994.

\bibitem{13}
\textsc{A.Meirmanov B.Zaltzman} \emph{Global in time solution to
the Hele-Shaw problem with a change of topology.} EJAM, Vol. 13,
pp. 431-447, 2002.

\bibitem{8}
\textsc{G.A.Omel'yanov} \emph{Dynamics and interaction of
nonlinear waves: multidimensional case.} Int. Conf.
"Differential Equations and Related Topics", dedicated to the
Centenary Anniversary of I.G.Petrovskii. Book of abstracts. Moscow
Univ. Press, 2001, pp. 305-306.

\bibitem{3}
\textsc{E.V.Radkevich} \emph{The Gibbs-Thomson correction and
conditions for the classical solution of the modified Stefan
problem.} 1991, Soviet Math. Doklady 43(1).

\bibitem{14}
\textsc{A.Friedman} \emph{Partial Differential Equations of
Parabolic Type.} Prentice-Hall, N.J., 1964.

\end{thebibliography}
\end{document}